\documentstyle[12pt]{article}
\newtheorem{theorem}{Theorem}[section]
\newtheorem{lemma}[theorem]{Lemma}

\newtheorem{proposition}[theorem]{Proposition}
\newtheorem{remark}[theorem]{Remark}
\newtheorem{definition}[theorem]{Definition}

\def\qed{\hspace{.5cm}$\dashv$}
\def\proof{\noindent{\em Proof. }}
\def\kappa{K_{S}^{2}}
\def\pro{I\!\!P}

\def\bbbc{{\mathchoice {\setbox0=\hbox{$\displaystyle\rm C$}\hbox{\hbox
 to0pt{\kern0.4\wd0\vrule height0.9\ht0\hss}\box0}}
 {\setbox0=\hbox{$\textstyle\rm C$}\hbox{\hbox
 to0pt{\kern0.4\wd0\vrule height0.9\ht0\hss}\box0}}
 {\setbox0=\hbox{$\scriptstyle\rm C$}\hbox{\hbox
 to0pt{\kern0.4\wd0\vrule height0.9\ht0\hss}\box0}}
 {\setbox0=\hbox{$\scriptscriptstyle\rm C$}\hbox{\hbox
 to0pt{\kern0.4\wd0\vrule height0.9\ht0\hss}\box0}}}} 

 \def\bbbz{{\mathchoice {\hbox{$\sf\textstyle Z\kern-0.4em Z$}}
 {\hbox{$\sf\textstyle Z\kern-0.4em Z$}}
 {\hbox{$\sf\scriptstyle Z\kern-0.3em Z$}}
 {\hbox{$\sf\scriptscriptstyle Z\kern-0.2em Z$}}}}
\newcommand{\nc}{\newcommand}
\nc{\cH}{{\cal H}}
\nc{\cA}{{\cal A}}
\nc{\cG}{{\cal G}}
\nc{\cC}{{\cal C}}
\nc{\cO}{{\cal O}}
\nc{\cI}{{\cal I}}
\nc{\cB}{{\cal B}}
\nc{\cY}{{\cal Y}}
\nc{\cK}{{\cal K}}
\nc{\cX}{{\cal X}}
\nc{\cS}{{\cal S}}
\nc{\cE}{{\cal E}}
\nc{\cF}{{\cal F}}
\nc{\cZ}{{\cal Z}}
\nc{\cQ}{{\cal Q}}
\nc{\cN}{{\cal N}}
\nc{\cP}{{\cal P}}
\nc{\cL}{{\cal L}}

\nc{\fr}{{\rightarrow}} 

\title{ABELIAN COVERS AND ISOTRIVIAL CANONICAL FIBRATIONS}
\author{Francesco Zucconi\\
 Dipartimento di Matematica e Informatica\\
 Universit\`{a} di Udine\thanks{Research carried out under the project "Geometria
Algebrica, Algebra Commutativa e aspetti comutazionali" (coordinatore nazionale
Claudio Pedrini).}}
\date{}

\begin{document}

\maketitle 
\begin{abstract}
We give a pure algebraic method to construct all the infinite families of surfaces
$S$ with isotrivial canonical fibration where $S$ is the minimal desingularization of
$X=Z/G$ and $G$
is an Abelian group acting diagonally on the product of two smooth curves: $Z=F\times
D$. In particular we recover all the known infinite families of surfaces  with
isotrivial canonical fibration and we produce many new ones. Our method
works in every dimension and, with minor modifications, it can be applied to
construct surfaces with canonical map of degree
$>1$.
\end{abstract}      
{\bf{Introduction.}}\par\noindent

An infinite family $\cF=\{S_{t}\}_{t\in\, T}$ of surfaces is a
set which intersects an infinite number of irreducible
components of the moduli space of all surfaces of general type: that is, for
every natural number $N\in\bbbz^{+}$ there exists $t\in T$ such that
$\chi(K_{S_{t}})>N$, where $K_{S_{t}}$ is the
canonical sheaf, $\chi(K_{S_{t}})=1-q(S_{t})+p_{g}(S_{t})$,
$p_{g}(S_{t})={\rm{dim}}_{\bbbc}H^{0}(S_{t},K_{S_{t}})$,
$q(S_{t})={\rm{dim}}_{\bbbc}\,H^{1}(S_{t},K_{S_{t}})$.

After Beauville \cite{Be}
showed that for some families of surfaces with $p_{g}\gg 0$ the
image $B$ of the canonical map
$\Phi_{\mid K_{S}\mid}:\,S-\rightarrow
I\!\!P(H^{0}(S,K_{S})^{\star})$ is a curve, it was 
worthwhile taking the trouble to establish a geography of these 
surfaces; other authors \cite{X1},\cite{X2}, \cite{K}, \cite{Xs},
\cite{MY}, \cite{Z1} developed the project. It was known soon that: $a)$ if 
$B$ is a curve and $\chi(K_{S})\geq 21$ then $\Phi_{\mid K_{S}\mid}$ is a
morphism and $2\leq g\leq 5$ where $g$ is the genus of the fibre, see
\cite{Be}[proposition 2.1]; 
$b)$ $B$ is the rational normal curve of
degree  
$p_{g}-1$ or the elliptic normal curve of degree  $p_{g}$ and $0\leq q(S)\leq
2$, \cite{X2}. So far some infinite families of surfaces with canonical
fibration 
$\phi_{\mid K_{S}\mid}:S\rightarrow B$ of genus $g=2$ or $g=3$ have been
found; they split into two cases: $i)$ $\phi_{\mid K_{S}\mid}:S\rightarrow
B$ is an isotrivial fibration, that is, the smooth fibres are
isomorphic each other; $ii)$ it is not an isotrivial
fibration; for this second class we refer to
\cite{Ca1} where the reader can find a complete bibliography.

The purpose of this article is to present
a uniform and algebraic procedure to construct infinite families of
surfaces with isotrivial canonical fibration; see Theorem \ref{otto}. 
Moreover we will show that the case
$g=5$ does not occur as Xiao conjectured in \cite{X3} and we will construct
explicitly all the genus $\leq 3$ isotrivial canonical fibrations $f:S\rightarrow B$
where
$S$ is the minimal desingularization of $X=Z/G$ and $G$ is an Abelian
group acting diagonally on the product surface $Z$. All the known examples of infinite
families can be recovered from our method and it can be generalized naturally to
varieties of dimension
$>2$. Finally we should point out the rich and unknown geometry that the following
two theorems show.\\ {\mbox{\bf{Theorem (The case $\mid G\mid\leq 8$):}}}\\ 
{\it{ Let $S$ be a minimal desingularization of the quotient $Z/G$ where $Z$ is the
product surface $F\times D$ and $G$ is an Abelian group acting faithfully on the
curves $F$, $D$ and diagonally on $Z$; set $A=F/G$ and $B=D/G$. Assume that the
canonical map of $S$ factorizes through $\pi_{2}:X\rightarrow B$. If
$p_{g}(S)=m-1\geq 3$, $\mid G\mid\leq 8$ and $g(F)\leq 3$ then only the following
cases occur:
\begin{equation}\label{zero}\begin{tabular}{|c|c|c|c|c|c|}\hline
  $g(A)$ & $g(B)$ &$G$  & $g(F)$& $g(D)$&
$K_{S}^{2}$\\
\hline

$0$ &$0$& ${\bbbz}/2\times{\bbbz}/2$ & $2$& $2m-1$& $4m$\\\hline

$0$ &$0$& ${\bbbz}/2\times{\bbbz}/2$ & $2$& $2m-2$& $4m-6$\\\hline

$0$ &$0$& ${\bbbz}/2\times{\bbbz}/2$ & $2$& $2m-3$& $4m-8$\\\hline

$0$ &$1$& ${\bbbz}/2\times{\bbbz}/2$ & $2$& $2m-1$& $4m-4$\\\hline

$1$ &$0$& ${\bbbz}/2$ & $2$& $m$& $4m-4$\\\hline

$0$ &$0$& ${\bbbz}/2\times{\bbbz}/2$ & $3$& $2m-1$& $8m$\\\hline

$0$ &$0$& ${\bbbz}/2\times{\bbbz}/2$ & $3$& $2m-2$& $8m-12$\\\hline

$0$ &$0$& ${\bbbz}/2\times{\bbbz}/2$ & $3$& $2m-3$& $8m-16$\\\hline

$0$ &$0$& ${\bbbz}/2\times{\bbbz}/2\times{\bbbz}/2$ & $3$& $4m+1$&
$8m$\\\hline

$0$ &$0$& ${\bbbz}/2\times{\bbbz}/2\times{\bbbz}/2$ & $3$& $4m-1$&
$8m-4$\\\hline

$0$ &$0$& ${\bbbz}/2\times{\bbbz}/2\times{\bbbz}/2$ & $3$& $4m-3$&
$8m-8$\\\hline

$0$ &$0$& ${\bbbz}/2\times{\bbbz}/2\times{\bbbz}/2$ & $3$& $4m-5$&
$8m-12$\\\hline

$0$ &$0$& ${\bbbz}/2\times{\bbbz}/2\times{\bbbz}/2$ & $3$& $4m-7$&
$8m-16$\\\hline

$0$ &$1$& ${\bbbz}/2\times{\bbbz}/2$ & $3$& $2m-1$& $8m-8$\\\hline

$0$ &$1$& ${\bbbz}/2\times{\bbbz}/2\times{\bbbz}/2$ & $3$& $4m-3$&
$8m-8$\\\hline

$1$ &$0$& ${\bbbz}/2\times{\bbbz}/2$ & $3$& $2m-1$& $8m-8$\\\hline

$1$ &$0$& ${\bbbz}/2\times{\bbbz}/2$ & $3$& $2m-2$& $8m-12$\\\hline

$1$ &$0$& ${\bbbz}/2\times{\bbbz}/2$ & $3$& $2m-3$& $8m-16$\\\hline

$2$ &$0$& ${\bbbz}/2$ & $3$& $m-1$& $8m-16$\\\hline
\end{tabular}\end{equation}\\}}
\noindent
If $\mid G\mid> 8$ only the case
$G={\bbbz}/2\times{\bbbz}/8$ occurs, but it produces many different
infinite families; see the proof of \ref{eccolotto} and the proof of
\ref{eccolottouno} for the last entry of the next table. More precisely:\\|\\\\
 
\noindent{\mbox{\bf{Theorem (The case $\mid G\mid> 8$):}}}\\
{\it{If $\mid G\mid> 8$ then $G={\bbbz}/2\times{\bbbz}/8$,
$g(F)=3$, $g(A)=0$. Moreover the
actions are completely described and if $m\geq 4$ only the following invariants occur:

\begin{equation}\label{mostro}\begin{tabular}{|c|c|c|c|c|}\hline
$g(B)$ &$p_{g}$  & $g(D)$& $K_{S}^{2}$\\\hline

$0$ &$m$ &  $8m-1$& $8m-6$\\\hline

$0$ &$m-1$ &  $8m-5$& $8m-16$\\\hline

$0$ &$m$ & $8m-5$& $8m-10$\\\hline

$0$ &$m-1$ & $8m-9$& $8m-20$\\\hline

$0$ &$m+1$ & $8m+4$& $8m$\\\hline

$0$ &$m$ & $8m$& $8m-4$\\\hline

$0$ &$m-1$ & $8m-4$& $8m-8$\\\hline

$0$ &$m$ & $8m$& $8m-3$\\\hline

$0$ &$m-1$ & $8m-4$& $8m-7$\\\hline

$0$ &$m-2$ & $8m-8$& $8m-11$\\\hline

$0$ &$m$ & $8m+3$& $8m$\\\hline

$0$ &$m-1$ & $8m-1$& $8m-4$\\\hline
$0$ &$m$ & $8m-1$& $8m-4$\\\hline

$0$ &$m-1$ & $8m-5$& $8m-8$\\\hline

$0$ &$m-1$ & $8m-9$& $8m-12$\\\hline

$0$ &$m+1$ & $8m+12$& $8m+7$\\\hline

$0$ &$m$ & $8m+8$& $8m+3$\\\hline

$0$ &$m-1$ & $8m+4$& $8m-1$\\\hline

$1$ &$m$ & $8m+1$& $8m$\\\hline
\end{tabular}\end{equation}\\
\noindent}}\\
{\bf{Acknowledgment:}}
I would like to thank F.Catanese for pointing out some evidence \cite{Ca} which
brought me to theorem \ref{otto}, then M. Manetti for showing me a short way
to exclude the genus-$5$ case and the group of geometers of Dipartimento di
Matematica del Politecnico di Torino, where I wrote the first
version of this paper, for their supportive attitude.\\ 
{\bf{Notations.}}\\
$S$ will be a surface with canonical map composed with the pencil
$f:S\rightarrow B$ and also it will be the desingularization of
$X=F\times D/G$ where
$G$ is an Abelian group acting faithfully on the two curves $F$, $D$ with quotients
$F/G=A$ and $D/G=B$ and acting diagonally on $F\times D$. In particular $f$ will
factorize through the projection
$X\rightarrow B$ and $F$ will be its fibre. Since we will be concerned with product
surfaces we will denote by capital letters also the points of
$A$, $B$, $D$, $B$ whenever it can arise no confusion.




\section{The diagonal action}

An action of a group $G$ on a variety $Z$ is called faithful if no non-trivial element
of $G$ acts trivially on $Z$. Let us consider an analytic faithful action of a finite
group $G$ on $Z$. The isotropy subgroup (stabilizer) of a point
$P\in Z$ is the subgroup $G_{P}$ of elements in $G$ fixing
$P$. An action is free if $G_{P}=\{1\}$ for every $P\in Z$. Let $X=Z/G$ be the
quotient variety with the natural projection $\pi_{G}:Z\rightarrow X$. A point
$P\in Z$ is called a ramification point if $\pi_{G}$ is ramified at $P$ and
$\pi_{G}(P)=[P]$ is called a branch point. It is well-known that all points in 
$\pi_{G}^{-1}[P]$ have the same multiplicity by $\pi_{G}$ and this number is called
the branching order of $[P]$. If $G$ acts on
two smooth curves $F$, $D$, we can define the diagonal action on
$F\times D$ by $(x_{1},x_{2})\mapsto (\gamma x_{1},\gamma x_{2})$ for
all $(\gamma,x_{1},x_{2})\in G\times F\times D$. We consider a finite
group $G$ acting faithfully on two smooth curves $F$,
$D$ and we denote $Z=F\times
D$,
$A=F/G$, $B=D/G$. Set
$X=Z/G$ for the quotient by the diagonal action and denote by 
$\pi_{A}:X\rightarrow A$, $\pi_{B}:X\rightarrow B$ the two projections. In
particular $\pi_{B}:X\rightarrow B$ is a constant moduli fibration meaning that all
the smooth fibres are isomorphic to $F$. A similar description applies to
$\pi_{A}:X\rightarrow A$: its smooth fibres are isomorphic to $D$. Since the
stabilizer of a  point $P$ of a curve is a {\bf{cyclic}} subgroup of $G$ then
the finite map $\pi_{G}:Z\rightarrow X$ is branched at the isolated
points $[(x_{1},x_{2})]$ such that 
$[(\gamma x_{1},\gamma x_{2})]=[(x_{1},x_{2})]$ for some $\gamma\in
G\setminus\{{\rm{id}}\}$. These points are precisely the singular locus of $X$, in
view of the purity of branch locus theorem. In particular $X$ has cyclic quotient
singularities, hence Hirzebruch-Jung singularities only. Then if $\nu :S\rightarrow X$
is the minimal desingularization of $X$, [\cite{Ser} p. 64], the fibre of $\nu$
over each singular point of $X$ is an Hirzebruch-Jung string [cf.\cite{BPV}]. Such
desingularization is called a {\it{ standard isotrivial surface}} after Serrano
\cite{Ser} and the induced fibrations $f_{1}:S\rightarrow A$, $f_{2}:S\rightarrow B$
are called {\it{ standard isotrivial fibrations}}. In the rest of this paper $G$
will be an Abelian group acting faithfully on the smooth curves
$F$, $D$ and we will also assume $g(F)\geq 2$ and $g(D)\geq 2$.
Let $\Delta:G\times Z\rightarrow Z$ be the diagonal action with quotient $Z/G=X$ and
with quotient map
$\rho :Z\rightarrow X$. Since $G$ acts freely on $Z$
outside a finite set of points $\{P_{1},...P_{t}\}\subset Z$ with non-trivial
stabilizer $G_{P_{i}}$, $X$ is a normal surface with only isolated rational
singularities [\cite{Bri} Satz 1.7]. Let
${\rm{Sing}}(X)=\{\rho (P_{1}),...\rho (P_{t})\}\subset Z$ be the singular locus of
$X$ and let
$j:X_{0}\rightarrow X$ be the natural inclusion. By well-known results on isolated
quotient singularities cf.\cite{Pi}[Th\'{e}or\`{e}me pp. 169] we have
${\cal{O}}_{X}=j_{\star}{\cal{O}}_{X_{0}}=(\rho_{\star}{\cal{O}}_{Z})^{G}$ and
$\omega_{X}=j_{\star}\omega_{X_{0}}=(\rho_{\star}\omega_{Z})^{G}$. Moreover, since $G$
is Abelian, there is an action $\alpha:G\times X\rightarrow X$
induced by $\alpha^{'}:G\times Z\rightarrow Z$ defined by
${\alpha^{'}}:(g;(x,y))\mapsto (gx,y)$.  We call $\alpha$ {\bf{the quotient action}}.
We set $Y=X/G$ and by a trivial computation we obtain $Y\, =\,A\times B$.

The next proposition is an extension of proposition 4.1 in \cite{P} to a surface
with isolated rational singularities.
\begin{proposition}\label{tre}
Let $\alpha: G\times X\rightarrow X$ be an action of an Abelian finite group $G$ on a
surface $X$ with isolated rational singularities and let
$\pi_{G}:X\rightarrow X/G$ be the projection. Set
$X/G=Y$ and assume that $Y$ is a smooth surface. 
Then for all $\chi\in G^{\star}$, where $G^{\star}$ is the group of
characters, there exists an invertible sheaf $L_{\chi}$ on
$Y$ such that

$$
\begin{array}{cc}
i)& \pi_{G_{\star}}{\cal{O}}_{X}=\oplus_{\chi\in G^{\star}}L^{-1}_{\chi}\\
\end{array}
$$
where $G$ acts on $L^{-1}_{\chi}$ via the character $\chi$ and the invariant summand 
$L_{id}={\cal{O}}_{Y}$. Furthermore we also have a decomposition into invertible
sheaves of $\pi_{G_{\star}}\omega_{X}$ such that:

$$
\begin{array}{cl}
ii)& (\pi_{G_{\star}}\omega_{X})^{\rm{inv}}=\omega_{Y}\\
iii) & (\pi_{G_{\star}}\omega_{X})^{\chi}=\omega_{Y}\otimes L_{\chi^{-1}}\\
\end{array}
$$
where $G$ acts on $(\pi_{G_{\star}}\omega_{X})^{\chi}$ via the character $\chi$.
\end{proposition}
\proof i). Since $X$ is normal, $Y$ is smooth and $G$ is Abelian then by \cite{P},
$\pi_{G}$ is flat and the action $\alpha$ induces the desired splitting on
$\pi_{G_{\star}}{\cal{O}}_{X}$.

ii). Let ${\rm{Sing}}(X)$ be the union of the singular points of $X$, 
$X_{0}=X\setminus{\rm{Sing}}(X)$, $Y_{0}=X_{0}/G$ and let $j:X_{0}\rightarrow X$,
$i:Y_{0}\rightarrow Y$,  $\pi_{G}^{0}:X_{0}\rightarrow Y_{0}$ be  respectively the
natural inclusions and the projection. Since $\pi_{G}^{0}$ is a smooth
$G$-cover then by \cite{P}[prop 4.1] there exist $L^{0}_{\chi}\in{\rm{Pic}}(Y_{0})$
such that $\pi^{0}_{G_{\star}}{\cal{O}}_{X_{0}}=\oplus_{\chi\in
G^{\star}}(L^{0}_{\chi})^{-1}$ and
$(\pi^{0}_{G_{\star}}\omega_{X_{0}})^{\chi}=\omega_{Y_{0}}\otimes L_{\chi^{-1}}^{0}$.
Since $Y$, $Y_{0}$ are smooth and $Y\setminus Y_{0}$ is a finite set of points then
by part i) we have $i^{\star}L_{\chi}=L^{0}_{\chi}$ and 
$i^{\star}\omega_{Y}=\omega_{Y_{0}}$
then 
$i_{\star}(\omega_{Y^{0}}\otimes L^{0}_{\chi^{-1}})=\omega_{Y}\otimes L_{\chi^{-1}}$.
We want to show that:
$$
(\pi_{G_{\star}}\omega_{X})^{\chi}=\omega_{Y}\otimes L_{\chi^{-1}}.
$$
Since $\omega_{X}=j_{\star}\omega_{X_{0}}$ then 
$i_{\star}\pi^{0}_{G_{\star}}\omega_{X_{0}}=\pi_{G_{\star}}\omega_{X}$.
On the other hand
$i_{\star}((\pi^{0}_{G_{\star}}\omega_{X_{0}})^{\chi})$
and $(i_{\star}\pi^{0}_{G_{\star}}\omega_{X_{0}})^{\chi}$ coincide outside a finite
set of points and $i_{\star}((\pi^{0}_{G_{\star}}\omega_{X_{0}})^{\chi})=
(i_{\star}\pi^{0}_{G_{\star}}\omega_{X_{0}})^{\chi}$ since $Y$ is smooth (normal).
Then $(\pi_{G_{\star}}\omega_{X})^{\chi}=i_{\star}((\omega_{Y^{0}}\otimes
L^{0}_{\chi^{-1}})=\omega_{Y}\otimes
L_{\chi^{-1}}$.\qed\\
\\
If $G$ is an Abelian finite group acting faithfully on two smooth curves $F$
and $D$, the diagonal $G$-action on the product surface $Z=F\times
D$ and the quotient action on the surface $X=Z/G$ forces some
structure on the invariants of $X$ which can be read from the $G$-action on 
$F$ and $D$. Then we are led naturally to the following definition:

\begin{definition}{\label{definizione chiave}}Let $G$ be an Abelian group acting
faithfully on two smooth curves, $F$ and $D$ where $g(F)\geq
2$, $g(D)\geq 2$. Let $Z=F\times D$ and let $\Delta:G\times
Z\rightarrow Z$ be the diagonal action with quotient $Z/G=X$. We will say that $X$ is
a {\bf{$G$- sandwich surface}} with top $F\times D$ and base
$A\times B$, where $A=F/G$, $B=D/G$. 
\end{definition}
The next
lemma relates 
the invariants of a $G$- sandwich surface $X$ with the invariants of
the curves which produce $X$.
\begin{lemma}\label{sandwich}
If $X$ is a $G$- sandwich surface with top $Z=F\times
D$ and base $A\times B$ then

$$\begin{array}{rl}
H^{0}(X,\Omega^{1}_{X})= & H^{0}(A,\Omega^{1}_{A})\oplus H^{0}(B,\Omega^{1}_{B})\\
\\
H^{0}(X,\omega_{X})= & H^{0}(Z,\omega_{Z})^{G}\\
\end{array}$$

\end{lemma}
\proof Let $\nu:W\rightarrow X$ be a minimal resolution of $X$. By \cite{Fre}[satz
1, p.99] we have 
$$\begin{array}{rl}
H^{0}(W,\Omega^{1}_{W})= & H^{0}(A,\Omega^{1}_{A})\oplus
H^{0}(B,\Omega^{1}_{B});\\
\\
H^{0}(W,\omega_{W})= & H^{0}(Z,\omega_{Z})^{G}.\\
\end{array}$$
Since $X$ has rational singularities then
$\nu_{\star}\Omega^{1}_{W}=i_{\star}\Omega^{1}_{X_{0}}=_{{\rm{def}}}\Omega^{1}_{X}$
and
$\nu_{\star}\omega_{W}=\omega_{X}$. The claim is now obvious.\qed\\\\
If $X$ is a $G$-sandwich surface with top $Z=F\times D$
and $H^{0}(F,\Omega^{1}_{F})=\oplus_{\chi\in G^{\star}}V_{1,\chi}$,
$H^{0}(D,\Omega^{1}_{D})=\oplus_{\chi\in G^{\star}}V_{2,\chi}$ are
the decompositions in subspaces where $G$ acts on
$V_{i,\chi}$ via the character $\chi$, there is a nice decomposition on the vector
space $H^{0}(X,\omega_{X})$ which is a sort of K\"{u}nneth formula for an Abelian
quotient of a direct product.
\begin{theorem}\label{quattro} 
Let $X=Z/G$ be a $G$-sandwich surface with top $Z=F\times D$ 
and base $Y=A\times B$ and let 
$H^{0}(F,\Omega^{1}_{F})=\oplus_{\chi\in G^{\star}}V_{1,\chi}$, 
$H^{0}(D,\Omega^{1}_{D})=\oplus_{\chi\in G^{\star}}V_{2,\chi}$ be the
decompositions in subspaces where $G$ acts on $V_{i,\chi}$ via the character $\chi$.
Then 
$$ H^{0}(X,\omega_{X})=\oplus_{\chi\in G^{\star}} V_{1,\chi}\otimes V_{2,\chi^{-1}}.$$
\end{theorem}
\proof
By \ref{sandwich} $H^{0}(X,\omega_{X})=H^{0}(Z,\Omega^{2}_{Z})^{G}$. On the other
hand by the K\"{u}nneth
formula $H^{0}(Z,\Omega^{2}_{Z})=(H^{0}(A,\Omega^{1}_{A})\otimes
H^{0}(B,\Omega^{1}_{B}))$, then $H^{0}(Z,\Omega^{2}_{Z})^{G}=(H^{0}(A,\Omega^{1}_{A})\otimes
H^{0}(B,\Omega^{1}_{B}))^{G}$, that is 
$H^{0}(X,\omega_{X})=
(\oplus_{\chi,\chi^{'}\in G^{\star}}V_{1\chi}\otimes V_{2\chi^{'}})^{G}=
\oplus_{\chi\in G^{\star}}V_{1\chi}\otimes V_{2\chi^{-1}}.$ 
                                    \qed\\\\
We have looked at the roof of X. Now we pay attention to the floor and we discover
that the pieces of the decomposition
$H^{0}(X,\omega_{X})=\oplus_{\chi\in G^{\star}} V_{1,\chi}\otimes V_{2,\chi^{-1}}$
induced by the diagonal action can be related to the pieces on $Y$ produced by
the quotient action $\alpha$.

\begin{proposition}\label{cinque} With the same hypotheses as for Theorem
\ref{quattro}, let
$\pi_{G}:X\rightarrow Y$ be the projection obtained by the quotient action $\alpha$
with associated splitting:
$\pi_{G_{\star}}{\cal{O}}_{X}=\oplus_{\chi\in G^{\star}}L^{-1}_{\chi}$. Then for all
$\chi\in G^{\star}$ we have
$$ H^{0}(X,\omega_{X})^{\chi}\simeq H^{0}(Y, \omega_{Y}\otimes
L_{\chi^{-1}})\simeq V_{1,\chi}\otimes V_{2,\chi^{-1}}.$$
\end{proposition}

\proof By \ref{tre} 
$(\pi_{G_{\star}}\omega_{X})^{\chi}=\omega_{Y}\otimes L_{\chi^{-1}}$. Then
$H^{0}(X,\omega_{X})^{\chi}\simeq H^{0}(Y, \omega_{Y}\otimes L_{\chi^{-1}})$. On the
other hand the action $\alpha$ is induced by the action $\alpha^{'}:G\times
Z\rightarrow Z$ which operates trivially on the second factor. Then by \ref{quattro}
we have  $H^{0}(X,\omega_{X})^{\chi}=V_{1,\chi}\otimes V_{2,\chi^{-1}}$; that is $
H^{0}(Y, \omega_{Y}\otimes L_{\chi^{-1}})=V_{1,\chi}\otimes V_{2,\chi^{-1}}$.\qed


\section{Surfaces with canonical map composed with a pencil}

Surfaces do
not usually come equipped with a fibration, but when they do and if the
fibration is the map associated with the canonical linear system, the
interplay between the genus and other additional structure on the fibres,
the genus of the base curve and the invariants naturally attached to the
surface provide a great deal of information.

A surjective morphism with connected fibre $F$, $f:S\rightarrow B$ of a
smooth projective surface $S$ onto a genus $b$ smooth curve $B$ is called a genus
$b$ pencil of curves of genus $g$ if $g$ is the arithmetic genus of $F$. From
now on we assume that $S$ is a surface of general type with
$p_{g}\geq 2$. We say that the canonical linear system $\mid K_{S}\mid$ is composed
with a pencil if the canonical image $\Sigma=\Phi_{\mid K_{S}\mid}(S)$ is a curve; in
this case let $\epsilon :S^{'}\rightarrow S$ be the elimination of the base points of
the moving part of $\mid K_{S}\mid$, then, taking the Stein factorization, we get
a genus-$b$ pencil $f: S^{'}\rightarrow B$ and a morphism $B\rightarrow \Sigma$. We
call the induced fibration $f: S^{'}\rightarrow B$ the {\bf{canonical fibration}}.
Since we will deal with singular surfaces we will recall Beauville's
definition of canonical map for a singular surface of general type,
\cite{Be}[pp.127]:

\begin{definition}\label{sette}
Let $X$ be a singular surface. Let $\nu :S\rightarrow X$ be a birational morphism 
where $S$ is a smooth surface. The rational map $\Phi_{\mid
K_{X}\mid}=_{\rm{def}}\Phi_{\mid K_{S}\mid}\circ (\nu)^{-1}$ is called the canonical
map of $X$. Moreover we put $p_{g}(X)=p_{g}(S)$ and one says that $X$ is of general
type if $S$ is of general type.
\end{definition}
Let us assume that the surface $X$ is of general type; by the unicity of the minimal
model, the definition  of $\Phi_{\mid K_{X}\mid}$ does not depend
on the choice of $\nu$. In particular $\Phi_{\mid K_{X}\mid}(X)$ is a curve if and
only if $\Phi_{\mid K_{S}\mid}(S)$  is a curve. Now we define a suitable class of
surfaces with canonical map composed with a pencil. This class includes all the known
examples of surfaces with isotrivial canonical fibration;  
see: \cite{Be}, \cite{X1}, \cite{X2}, \cite{K}, \cite{Z2}. 

\begin{definition}\label{g-canonical}
We shall say that $S$ is a $G$-sandwich canonically $g$-fibred surface on $B$ with
top $F\times D$ and base $A\times B$ ($G$-
canonical for short) if 1) $S$ is a smooth model of a $G$- sandwich surface $X$ with
top
$F\times D$ and base $A\times B$ [see \ref{definizione
chiave}]; 2) $\Phi_{\mid K_{S}\mid}(S)=B$, $g=g(F)$; 3) the induced
fibration $f:S\rightarrow B$ factorizes through $\pi_{B}:X\rightarrow B$.
\end{definition}
In the following theorem we will prove an effective method to characterize
$G$-sandwich canonically $g$-fibred surfaces. This theorem relies on theorem
\ref{quattro}.

\begin{theorem}\label{otto}
Let $X$ be a $G$-sandwich surface with top $Z=F\times D$ and base 
$Y=A\times B$ with $p_{g}\geq 2$ and let $\nu:S\rightarrow X$ be a minimal
desingularization of $X$. Let
$H^{0}(F,\Omega^{1}_{F})=\oplus_{\chi\in
G^{\star}}V_{1,\chi}$, $H^{0}(D,\Omega^{1}_{D})=\oplus_{\chi\in
G^{\star}}V_{2,\chi}$ be the decompositions into subspaces where $G$ acts on
$V_{i,\chi}$, $i=1,2$, via the character $\chi$. Then $S$ is a $G$-sandwich
canonically
$g$-fibred surface on
$B$ if and only if there exists an unique $\chi_{0}\in G^{\star}$ such that the
following two conditions hold:
$$
\begin{array}{cl}
i) & V_{1,\chi_{0}}\otimes V_{2,\chi_{0}^{-1}}\neq 0\\ ii) &
{\rm{dim}}_{\bbbc}V_{1,\chi_{0}}=1,  {\rm{dim}}_{\bbbc}V_{2,\chi^{-1}_{0}}=p_{g}.\\
\end{array}
$$
\end{theorem}
\proof If there exists an unique $\chi_{0}\in G^{\star}$ such that 
${\rm{dim}}_{\bbbc}V_{1,\chi_{0}}\neq 0$,${\rm{dim}}_{\bbbc}V_{2,\chi^{-1}_{0}}\neq 0$
then by
\ref{tre} and \ref{cinque} $H^{0}(x,\omega_{X})^{\chi}=0$ for all $\chi\neq\chi_{0}$.
In particular by \ref{cinque} $H^{0}(X,\omega_{X})= V_{1,\chi_{0}}\otimes
V_{2,\chi_{0}^{-1}}$. Now, if ${\rm{dim}}_{\bbbc}V_{1,\chi_{0}}=1$ and $\phi ,
\psi\in H^{0}(X,\omega_{X})$, then $\psi /\phi=\mu (y)$ is a meromorphic function
which depends only on $y\in B$; whence $\Phi_{\mid K_{X}\mid}$ factors through
$f_{2}: X\rightarrow B$. Conversely, if $\Phi_{\mid K_{X}\mid}$  is composed with
the pencil $f_{2}: X\rightarrow B$ and there exist $\chi,\chi^{'}\in G^{\star}$
such that $V_{1,\chi}\otimes V_{2,\chi^{-1}}\neq 0$ and $V_{1,\chi^{'}}\otimes
V_{2,\chi^{'^{-1}}}\neq 0$ or there exists $\chi_{0}\in G^{\star}$ such that
$({\rm{dim}}_{\bbbc}V_{1,\chi_{0}})\times ({\rm{dim}}_{\bbbc}V_{2,\chi^{-1}_{0}})\geq
2({\rm{dim}}_{\bbbc}V_{2,\chi^{-1}_{0}})$ then there exist $\eta_{1}, \eta_{2}\in
H^{0}(F,\Omega^{1}_{F})$ and  $\sigma_{1}, \sigma_{2}\in
H^{0}(D,\Omega^{1}_{D})$ such that
$\omega_{i}=\eta_{i}\wedge\sigma_{i}$, $i=1,2$ are $2$-forms on $X$ and they induce a
moving linear system on the fibres $F$ of $f_{2}: X\rightarrow B$. In
particular
$\Phi_{\mid K_{X}\mid}(F)$ can not be a point; a contradiction. The same argument
shows that ${\rm{dim}}_{\bbbc}V_{1,\chi_{0}}=1$ and
${\rm{dim}}_{\bbbc}V_{2,\chi^{-1}_{0}}=p_{g}$.\qed
\\ 

Surfaces with canonical map composed with a 
pencil are very exceptional in the theory of surfaces of general type. In particular
there are strong bounds on some of their invariants. The next proposition is a
specialization of these bounds to the case of a $G$-sandwich canonically $g$-fibred
surface.

\begin{theorem}\label{novedieci}
If $S$ is a $G$-sandwich canonically $g$-fibred surface on $B$ with top
$F\times D$ and base $A\times B$ with $p_{g}(X)\geq 11$ and
$g=g(F)$ then $2\leq g\leq 5$ and only the following cases can occur:
$$A)\,\,
g(A)=q(X)\leq 2 \,\, {\rm{and}}\,\,  g(B)=0;$$
\noindent
or
$$B)\,\, g(A)=0 \,\, {\rm{and}}\,\,  g(B)=q(X)=1.$$
Moreover if $p_{g}\geq{\rm{max}}\{ 85-2q(X), 44\}$ then $g=5$ does not occur. 
\end{theorem}
\proof
By \ref{sette} there is a birational morphism $\eta :S\rightarrow X$ where $\Phi_{\mid
K_{S}\mid}$ is composed with the pencil $f=f_{2}\circ\eta$. It is easy to see that
the moving part of $\mid K_{S}\mid$ is without base points. In particular
$F$ is the generic fibre $F$ of $f$ and then the bound $2\leq g\leq
5$ on $g(F)$ is in \cite{Be}[Proposition 2.1]. Since $q(X)= g(A)+g(B)$
then $A)$ and $B)$ are easily obtained by \cite{X2}[proposition]. Now we show that
$g=5$ and $p_{g}\geq 44$ does not occur. In fact if $g=5$ then by
\cite{Xs} (see also \cite{Z2})
$K_{S}=8H+V+f^{\star}(\gamma)$ where
$p_{g}(S)=h^{0}(B,\cO_{B}(\gamma))$, $VF=0$, $HF=1$. In particular if
$g(B)=0$,
$K_{S}H\geq 3$ since $p_{g}\geq 44$ and $K_{S}H\geq 1$ if $g(B)=1$.
Then $K_{S}^{2}\geq 8(p_{g}(S)+2)$ if $g(B)=0$ and $K_{S}^{2}\geq 8(p_{g}(S)+1)$
if $g(B)=1$. On the other hand by
\cite{Ser}[Proposition 5.3] the Euler characteristic satisfies: $2c_{2}(S)\geq
K_{S}^{2}$, i.e. $c_{2}(S)\geq 4(p_{g}(S)+2)$ in case $A)$ and $c_{2}(S)\geq
4(p_{g}(S)+1)$ in case $B)$ : a contradiction with the Noether formula.\qed \\\\
>From Theorem \ref{novedieci} we see that to produce a $G$-sandwich canonically
$g$-fibred surface with large invariants we must try to construct surfaces with
$p_{g}(S)\gg 0$. On the other hand the above restriction on the genus
of $F$, $A$, $B$ enable us to build all the $G$-sandwich canonically
$g$-fibred surfaces. In fact the two conditions in theorem \ref{otto} show a
procedure to manifest all these surfaces. The following remark will be useful:

\begin{remark}\label{undici}If $S$ is a $G$-sandwich canonically
$g$-fibred surface on $B$ then there exist 
$\eta\in H^{0}(F,\Omega^{1}_{F})$, 
$\omega_{1},...,\omega_{p_{g}(S)}\in H^{0}(D,\Omega^{1}_{D})$ such
that 
$$\langle \eta\wedge\omega_{1},...,\eta\wedge\omega_{p_{g}(S)}\rangle$$
gives a basis of $H^{0}(S,\omega_{S})$. Moreover 
${\rm{div}}(\eta)\in {\rm{Div}}(F)$ is a $G$-invariant divisor on $F$.
\end{remark}

\proof The proof is trivial if in \ref{otto} we put $V_{1,\chi_{0}}=\eta\bbbc$ and 
$V_{2,\chi^{-1}_{0}}=\oplus_{i=1}^{\small{p_{g}(S)}}\omega_{i}\bbbc$.\qed
\\
\\
We will be concerned mostly with $G$-canonical surfaces where $G$ is a product of some
copies of $\bbbz/2$. In-fact we will use tacitly the following description of
quotient surfaces $X=Z/G$ where all the stabilizers are isomorphic to ${\bbbz}/2$.
\\ 

\begin{theorem}\label{barlow}
Let $Z$ be a nonsingular surface with $G$ acting freely outside the finite set of
points $\{P_{1},...P_{t}\}\subset Z$ with stabilizer $G_{P_{i}}={\bbbz}/2$ for
$i=1,...,t$. Let $\pi:Z\rightarrow X=Z/G$ be the quotient map and let $\nu:
S\rightarrow X$ be the minimal resolution (i.e. the blow up of the $t$ nodes of $X$).
Then if $n=\mid G\mid$ we have:
$$
\begin{array}{cl}
i) & K_{S}^{2}=K_{X}^{2}=\frac{1}{n}K_{Z}^{2}\\
ii) & \chi(\cO_{S})=\chi(\cO_{X})=\frac{1}{n}(\chi(\cO_{Z})+\frac{t}{4})\\
iii) & {\it{ If\, Z\, is\, a\,minimal\, surface\, of\, general\, type\, then\, so\,
is\, S.}}
\end{array}
$$
\end{theorem}
\proof See cf. \cite{Ba}[p. 295].\qed
\\
\\
Finally we recall that the aim of this paper is to look into {\bf{infinite families}}
of surfaces with isotrivial canonical map. By \cite{Be}[Proposition 1.7] the minimal
surfaces $S$ of general type with $\chi(\cO_{S})<N$ (or $K^{2}_{S}<N$ or
$p_{g}<N$) are in a limited family. Then, given $m\in\bbbz^{+}$, $m\gg 0$ we want
surfaces with $p_{g}=m$ and isotrivial canonical fibration. By \ref{novedieci} the
irregularity
$q=q(S)=g(A)+g(B)$ is at most $2$ and if $p_{g}\gg 0$, $2\leq
g=g(F)\leq 4$. All known examples are with $g=2$ or $g=3$.
We postpone the case $g=4$ and the exceptional limited families with $g\geq 5$ to a
forthcoming paper. We will classify all the
$G$-sandwich canonically
$g$-fibred surfaces with $p_{g}\gg 0$ and $g=2$ or $g=3$. The seven cases to classify
are listed below (see
\ref{novedieci}):

\begin{equation}\label{tabellacasi}
\begin{tabular}{|c|c|c|c|}\hline
{\rm{Case}} & $g(F)$ &  $g(A)$& $g(B)$\\ \hline
$B^{(2)}$ &$2$ & $0$ & $1$\\ \hline

$A^{(2)}_{0}$ &$2$ & $ 0$ & $ 0$\\ \hline

$A^{(2)}_{1}$ &$2$ & $ 1$ & $0$\\ \hline

$B^{(3)}$ &$3$ & $0$ & $1$\\ \hline

$A^{(3)}_{0}$&$3$ & $0$ & $0$\\ \hline

$A^{(3)}_{1}$ &$3$ & $1$ & $0$\\ \hline

$A^{(3)}_{2}$ &$3$ & $2$ & $0$\\ \hline
\end{tabular}
\end{equation}
\noindent


\section{$g(F)=2$}

The outline of the classification proof of infinite families of $G$-sandwich
canonically $2$-fibred surfaces is the same as that of genus-$3$ case. On the
other hand it is easier and it is clear how it works. We describe  all the
genus-2 infinite families even in the known cases. Moreover we will prove a unicity
theorem which was a gap in this theory.
We will need the following well-known result:
\begin{lemma}\label{unpuntoramificazione}
If $G$ is a finite Abelian group acting on a smooth curve $C$ then $C\rightarrow C/G$
has at least two branch points .
\end{lemma}
\proof
See cf. \cite{Ser}[Lemma 5.7].\qed
\\
\\
Theorem \ref{otto} works if we have classified all the $G$-actions on a curve
$F$ of genus $2$ with an invariant canonical divisor: see remark \ref{undici}.
In particular the splitting
$H^{0}(F,\Omega^{1}_{F})=\oplus_{\chi\in G^{\star}}V_{1,\chi}$ has to
be computed:

\begin{lemma}\label{tredici}
Let $G$ be an Abelian group acting on a smooth genus-2 curve $F$ with quotient
$F/G=A$. Let $\eta\in H^{0}(F,\Omega^{1}_{F})$ be a
$1$-form such that ${\rm{div}}(\eta)$ is a $G$-invariant divisor and let
$H^{0}(F,\Omega^{1}_{F})=\oplus_{\chi\in G^{\star}}V_{1,\chi}$ be the
decompositions into subspaces where $G$ acts on $V_{1,\chi}$ via the character
$\chi$. Then only the following cases occur: 
\begin{equation}\label{tabelladue}
\begin{tabular}{|l|c|l|}\hline
  $g(A)$ &  $G$  & $H^{0}(F,\Omega^{1}_{F})$\\ \hline
$1$ & ${\bbbz}/2$ & $V_{1,{\rm{id}}}\oplus V_{1,\chi}$\\ \hline

$0$ & ${\bbbz}/2$ & $V_{1,\chi}$\\ \hline

$0$ & ${\bbbz}/3$ & $V_{1,\chi}\oplus V_{1,\chi^{2}}$\\ \hline

$0$ & ${\bbbz}/4$ & $V_{1,\chi}\oplus V_{1,\chi^{3}}$\\ \hline

$0$ & ${\bbbz}/5$ & $V_{1,\chi}\oplus V_{1,\chi^{3}}$\\ \hline

$0$ & ${\bbbz}/6$ & $V_{1,\chi}\oplus V_{1,\chi^{5}}$\\ \hline

$0$ & ${\bbbz}/2\times {\bbbz}/2$ & $V_{1,\chi_{1}}\oplus V_{1,\chi_{2}}$\\ \hline 
\end{tabular}
\end{equation}
\noindent
where $\chi$ is a generator of $({\bbbz}/d)^{\star}$ if $d=2,3,4,5,6$ and  $\chi_{1}$, $\chi_{2}$ generate 
$({\bbbz}/2\times {\bbbz}/2)^{\star}$.
\end{lemma}
\proof The claim follows from \ref{unpuntoramificazione} and from the well-known
argument by Hurwitz to compute all the possible ramification indexes.\qed
\\
\\
Our aim is to classify all $G$-sandwich canonically $2$-fibred surfaces with
$p_{g}\gg 0$. We start with the case $B^{2}$ in table (\ref{tabellacasi}).

\begin{proposition}\label{quattordici}
If $S$ is a $G$-sandwich canonically $2$-fibred surface on an elliptic curve with 
$p_{g}(S)=m\geq 2$ then $S$ is the minimal desingularization of $X=Z/G$ where
$G={\bbbz}/2\times {\bbbz}/2$ and it acts diagonally on $Z=F\times D$,
$g(F)=2$, $g(D)=2m+1$, $g(A)=0$, $g(B)=1$. Moreover $X$ has
$8m$ nodes, $\Phi_{\mid K_{S}\mid}$ is composed with the pencil $S\rightarrow
B$ and $K^{2}_{S}=4\chi(\cO_{S})$.
\end{proposition}
\proof Let $S$ be a $G$-canonical surface with $\Phi_{\mid K_{S}\mid}$ composed 
with the pencil $f_{2}: S\rightarrow B$ where $g(B)=1$ and $g(F)=2$.
By \ref{undici} there is a $G$-invariant $1$-form on $F$. By
\ref{tredici} all these $G$ actions are known. Then we have to study all the
$G$-covers $D\rightarrow B$ such that the decomposition families 
$\{V_{1,\chi}\}_{\chi\in G^{\star}}$, $\{V_{2,\chi}\}_{\chi\in G^{\star}}$ satisfy
the condition of \ref{otto}. In particular we want the character $\chi_{0}$.
Since $g(B)=1$ then by \ref{novedieci}, $g(A)=0$. Now we follow the list
in
\ref{tredici} table (\ref{tabelladue}).

The case {\boldmath $G={\bbbz}/2$} is impossible. 
In fact  $H^{0}(F,\Omega^{1}_{F})=V_{1,\chi_{0}}$, then
${\rm{dim}}_{\bbbc}V_{1,\chi_{0}}=2$ contradicting \ref{otto}.

The cases {\boldmath $G={\bbbz}/d$} where $d=3,4,5,6$ are impossible. 
The proofs are the same. We explain the proof of $d=3$ only. The Galois map
$D\rightarrow B$ is given by
${\cal{L}}_{\chi},{\cal{L}}_{\chi^{2}}\in {\rm{Pic}}(B)$ and two effective
divisors $D_{1}$, $D_{2}$ such that:

\begin{equation}
\label{zetatre}
\left\{ \begin{array}{rl}
3{\cal{L}}_{\chi}= &D_{1}+2D_{2}\\
{\cal{L}}_{\chi^{2}} = & 2{\cal{L}}_{\chi}-D_{2}\\ 
\end{array}\right.
\end{equation}

\noindent
By \ref{otto} we have $V_{2,\chi}=0$ or $V_{2,\chi^{2}}=0$. 
Since $V_{2,\chi^{2}}=H^{0}(D,\Omega^{1}_{D})^{\chi^{2}}=
H^{0}(B,{\cal{O}}_{B}({\cal{L}}_{\chi}))\neq 0$ then $V_{2,\chi}=0$. Since
$V_{2,\chi}=H^{0}(B,{\cal{O}}_{B}({\cal{L}}_{\chi^{2}}))$, 
${\rm{deg}}({\cal{L}}_{\chi^{2}})=0$. We set 
${\rm{deg}}({\cal{L}}_{\chi})=l_{\chi}$, ${\rm{deg}}(D_{1})=d_{1}$,
${\rm{deg}}(D_{2})=d_{2}$. Taking the degrees, we have
$3l_{\chi}=d_{1}+2d_{2}$ and $2l_{\chi}=d_{2}$, that is $l_{\chi}=d_{1}=d_{2}=0$. In
particular $p_{g}(S)=0$.

{\boldmath $G={\bbbz}/2\times {\bbbz}/2$}. 
In this case the covering $D\rightarrow B$ is given by the
data ${\cal{L}}_{\chi_{1}},{\cal{L}}_{\chi_{2}}, {\cal{L}}_{\chi_{3}}\in
{\rm{Pic}}(B)$
 where $\{ 1,\chi_{1},\chi_{2},\chi_{3}\}=G^{\star}$ and three effective divisors on
$B$, $D_{1}$, $D_{2}$, $D_{3}$ such that:
  
\begin{equation}
\label{zetadueperse}
\left\{ \begin{array}{rl}
2{\cal{L}}_{\chi_{1}}= &D_{1}+D_{3}\\
2{\cal{L}}_{\chi_{2}} = & D_{2} +D_{3}\\
{\cal{L}}_{\chi_{3}}= & {\cal{L}}_{\chi_{1}}+{\cal{L}}_{\chi_{2}}-D_{3}\\ 
\label{eq:zetadueperse}
\end{array}\right.
\end{equation}

\noindent
By \ref{otto} and \ref{tredici}, $V_{2,\chi_{1}}=0$ or $V_{2,\chi_{2}}=0$. Assume
$0=V_{2,\chi_{1}}=H^{0}(B,{\cal{O}}_{B}({\cal{L}}_{\chi_{1}})$. In particular 
${\rm{deg}}({\cal{L}}_{\chi_{1}})=0$. Then, taking the degrees, we obtain 
$l_{\chi_{1}}=d_{1}=d_{3}=0$, $2l_{\chi_{2}}=d_{2}$, $l_{\chi_{2}}=l_{\chi_{3}}$.
Since $D_{i}$ is effective we have $D_{1}=D_{3}=0$ and ${\cal{L}}_{\chi_{1}}$ is a
$2$-torsion non-trivial divisor since $D$ is connected. Now if we set
$l_{\chi_{2}}=m$ we have the claim.\qed
 

\begin{remark} This unique family is described by Beauville in 
\cite{Be}[Example 2.6 p. 127].
\end{remark}
Now we describe the case $A_{1}^{2}$ in the table (\ref{tabellacasi}).
\begin{proposition}\label{sedici}
If $S$ is a $G$-sandwich canonically $2$-fibred surface on a rational curve with
$q(S)=1$, $p_{g}(S)=m\geq 2$ then $S$ is the minimal desingularization of $X=Z/G$
where $G={\bbbz}/2$ acts diagonally on $Z=F\times D$,
$g(F)=2$, $g(D)=m+1$, $g(A)=1$, $g(B)=0$. Moreover
$X$ has $4m$ nodes, $\Phi_{\mid K_{S}\mid}$ is composed with $S\rightarrow
B$ and $K^{2}_{S}=4\chi(\cO_{S})$.
\end{proposition}

\proof We suppose that $\mid K_{X}\mid$ is composed with 
$f_{2}: X\rightarrow B$ where $g(B)=0$. Since 
$1=q(X)=h^{0}(Z,\Omega^{1}_{Z})^{G}$ then $g(A)=1$. Furthermore, by definition,
$g(F)=2$ and therefore by \ref{tredici}, $G={\bbbz}/2$ and  
$H^{0}(F,\Omega^{1}_{F})=H^{0}(F,\Omega^{1}_{F})^
{\rm{id}}\oplus H^{0}(F,\Omega^{1}_{F})^{\chi}$ where $G^{\star}=\{1,
\chi\}$. From \ref{otto} it follows that each hyperelliptic covering
$D\rightarrow B$ branched on $2m+4$ points gives a solution. 
\qed


\begin{remark} This infinite family is described also by Beauville 
in \cite{Be}[Variante p. 126].
\end{remark}
We finish this section with the case $A^{2}_{0}$ in the table (\ref{tabellacasi}).
As far as we know there are examples for the case $A^{2}_{0}$ only in \cite{X1}.
\begin{proposition}\label{diciotto}
If $S$ is a $G$-sandwich canonically $2$-fibred surface on a rational curve, with
$q(S)=0$, $p_{g}(S)=m-1\geq 2$ then $S$ is the minimal desingularization of $X=Z/G$
where $G={\bbbz}/2\times{\bbbz}/2$ acts diagonally on $Z=F\times
D$, $g(F)=2$, $X$ has $t$ nodes and $S$ is in one of the following
three classes:
$$
\begin{array}{rcll} 
I) & g(D)=2m-1& K^{2}_{S}=4\chi(\cO_{S})& t=8m+8\\
II) & g(D)=2m-2 & K^{2}_{S}=4\chi(\cO_{S})-6& t=8m+12\\
III) & g(D)=2m-3& K^{2}_{S}=4\chi(\cO_{S})-8& t=8m+16.\\
\end{array}
$$
\end{proposition}
\proof
We can exclude the cases $G={\bbbz}/d$, $d=2,3,4,5,6$ of \ref{tredici} as 
in \ref{quattordici}. Let ${\boldmath{G={\bbbz}/2\times {\bbbz}/2}}$. Since
$B=\pro^{1}$ then 
${\rm{Pic}}(\pro^{1})=\bbbz$ and to obtain the solutions we need to solve the
degree system associated to (\ref{zetadueperse}) with the condition
$H^{0}(B,\omega_{B}({\cal{L}}_{\chi_{1}}))=0$. Since  
$H^{0}(B,\omega_{B}({\cal{L}}_{\chi_{1}}))\simeq
H^{0}(\pro^{1},{\cal{O}}_{\pro^{1}}(l_{\chi_{1}}-2))=0$ and
${\rm{deg}}({\cal{L}}_{\chi_{1}})>0$ the constraint is 
${\rm{deg}}({\cal{L}}_{\chi_{1}})=1$. Taking the degrees in (\ref{zetadueperse}) we
obtain: 
$$
\begin{array}{rl}
I) & d_{1}=2,
d_{2}=2m, d_{3}=0, l_{\chi_{1}}=1, l_{\chi_{2}}=m, l_{\chi_{3}}=m+1\\
II) & d_{1}=1,d_{2}=2m-1, d_{3}=1, l_{\chi_{1}}=1, l_{\chi_{2}}=m
l_{\chi_{3}}=m\\ III)& d_{1}=0,
d_{2}=2m-2, d_{3}=2, l_{\chi_{1}}=1, l_{\chi_{2}}=m, l_{\chi_{3}}=m-1,
\end{array}
$$
which are the desired solutions.\qed \\
We have  shown the cases with $g(F)=2$ in the table (\ref{zero}).

\subsection{The geometrical construction of the three families of \ref{diciotto}} We
want to construct the surfaces of \ref{diciotto}\label{stanco} in a more geometrical
way. We consider five distinct points
$Q_{1},Q_{2},Q_{3},Q_{4},Q_{5}\in I\!\!P^{1}=A$. Let
$\rho_{E}:E\rightarrow A$ and $\rho_{E^{'}}:E^{'}\rightarrow A$ be the
$\bbbz/2$- covers branched on 
$Q_{1},Q_{2},Q_{3},Q_{4}$ and on $Q_{2},Q_{3},Q_{4},Q_{5}$ respectively. Set
$\rho_{E}^{-1}(Q_{i})=Q_{i}^{E}$,
$i=1,2,3,4$, $\rho_{E}^{-1}(Q_{5})=P_{1}^{E}+P_{2}^{E}$ and
set $\rho_{E^{'}}^{-1}(Q_{i})=Q_{i}^{E^{'}}$,
$i=2,3,4,5$, $\rho_{E^{'}}^{-1}(Q_{1})=P_{1}^{E^{'}}+P_{2}^{E^{'}}$. Moreover let
$i_{E}:E\rightarrow E$, $i_{E^{'}}:E^{'}\rightarrow E^{'}$ be the involutions
associated to $\rho_{E}$, $\rho_{E^{'}}$ respectively. Let
$C$ be the normalization of the fibre product $E\times_{A}E^{'}$, thus
$g(C)=2$ and the natural map $\rho_{1}:C\rightarrow A$ is
a $G=\bbbz/2\times\bbbz/2$- cover. We call
$\rho_{1}^{-1}(Q_{i})=Q_{i1}^{C}+Q_{i2}^{C}$, $i=1,2,3,4,5$; in particular
$\rho_{1}^{\star}(Q_{i})=2Q_{i1}^{C}+2Q_{i2}^{C}$ where $i=1,2,3,4,5$. By abuse of
notation we put $G=\langle i_{E}, i_{E^{'}}\rangle$. Then $E=C/\langle
i_{E^{'}}\rangle$ and $E^{'}=C/\langle i_{E}\rangle$. It is easy to see that
$Q_{11}^{C}, Q_{12}^{C}$ are the ramification points of the natural map
$\pi_{E^{'}}:C\rightarrow E^{'}$ and $Q_{51}^{C}, Q_{52}^{C}$ are the ramification
points of $\pi_{E}:C\rightarrow E$. In particular $i_{E}(Q_{1j}^{C})=Q_{1j}^{C}$,
$j=1,2$; $i_{E}(Q_{i1}^{C})=Q_{i2}^{C}$,
$i=2,3,4,5$; $i_{E^{'}}(Q_{5j}^{C})=Q_{5j}^{C}$,
$j=1,2$; $i_{E^{'}}(Q_{i1}^{C})=Q_{i2}^{C}$, $i=1,2,3,4$. Then the six
points $Q_{ij}^{C}$ where $i=2,3,4$ and $j=1,2$ are fixed by the involution 
$i_{E}{i_{E^{'}}}$ and $\pi_{\Gamma}:C\rightarrow \Gamma=C/\langle
i_{E}i_{E^{'}}\rangle$ is the hyperelliptic map. Now we can find a basis 
$\langle \eta_{E}, \eta_{E^{'}}\rangle$ of $H^{0}(C,\Omega^{1}_{C})$
where $\eta_{E}$ is a $i_{E^{'}}$-invariant form and $\eta_{E^{'}}$
is $i_{E}$-invariant. If $G^{\star}=\{{\rm{id}},\chi_{1},\chi_{2},\chi_{3}\}$ where 
$\chi_{1}(i_{E})=-1$, $\chi_{1}(i_{E^{'}})=1$, $\chi_{2}(i_{E})=1$,
$\chi_{2}(i_{E^{'}})=-1$, $\chi_{1}\chi_{2}=\chi_{3}$ then $\langle
\eta_{E}\rangle=V_{1,\chi_{1}}$, $\langle\eta_{E^{'}}\rangle=V_{1,\chi_{2}}$. We have
construct the $\bbbz/2\times\bbbz/2$-action on the genus-$2$ curve $C=F$. We
will construct the $\bbbz/2\times\bbbz/2$-action on $D$ obtained
in \ref{diciotto}. We consider
${\widetilde{E}}\rightarrow I\!\!P^{1}=B$ and ${\widetilde{E^{'}}}\rightarrow
I\!\!P^{1}=B$ two $2$-to-$1$ morphisms branched on $P_{1},...P_{2m}$ and on
$A_{1},A_{2}$ respectively and let
$i_{{\widetilde{E}}}:{\widetilde{E}}\rightarrow{\widetilde{E}}$,
$i_{{\widetilde{E}}^{'}}:{\widetilde{E}}^{'}\rightarrow{\widetilde{E}}^{'}$ be the
natural involutions. Denote by $D$ the normalization
of ${\widetilde{E}}\times_{B} {\widetilde{E}}^{'}$; the
natural map $\rho_{2}:D\rightarrow B$ is a
$G=\bbbz/2\times\bbbz/2$- cover and  
$H^{0}(D,\Omega^{1}_{D})=V_{1,\chi_{2}}\oplus V_{1,\chi_{3}}$. In
fact it is easy to find a basis  $\langle\alpha_{1},...,\alpha_{2m-1}\rangle$ such
that $\langle
\alpha_{1},...,\alpha_{m-1}\rangle=
V_{1,\chi_{2}}=H^{0}({\widetilde{E}},\Omega^{1}_{{\widetilde{E}}})$, and 
$\langle
\alpha_{m},...,\alpha_{2m-1}\rangle=
V_{2,\chi_{3}}=H^{0}({\widetilde{\Gamma}},\Omega^{1}_{{\widetilde{\Gamma}}})$ where 
${\widetilde{\Gamma}}=D/\langle i_{{\widetilde{E}}}
i_{{\widetilde{E}}^{'}}\rangle$. The group $G=\langle
(i_{E},i_{{\widetilde{E}}}),(i_{E^{'}},i_{{\widetilde{E}}^{'}})\rangle$ acts
diagonally on the product surface $Z=F\times D$ and if
$X=Z/G$,
$H^{0}(X,\Omega^{1}_{X})=H^{0}(Z,\Omega^{1}_{Z})^{G}=0$,
$H^{0}(X,\Omega^{2}_{X})=H^{0}(Z,\Omega^{2}_{Z})^{G}= V_{1,\chi_{2}}\otimes
V_{2,\chi_{2}}$; that is
$H^{0}(X,\Omega^{2}_{X})\simeq
\langle\eta\wedge\alpha_{1},...\eta\wedge\alpha_{m-1}\rangle$. Then $\Phi_{\mid
K_{X}\mid}$ yields a pencil with image the rational normal curve of degree $m-1$ and
general fibre $F$. We remark that if $A_{i}\neq R_{j}$, $i=1,2$, $j=1,...,2m$
we obtain the family $I)$; if  $A_{1}=R_{1}$  and  $A_{2}\neq R_{j}$ $j=1,...,2m$, we
have $II)$; finally, if  $A_{1}=R_{1}$, $A_{2}=R_{2}$ then we have $III)$. We can
carry on the explicit computation of the nodes occurring in each
family. We recall that 
$\rho_{1}^{-1}(Q_{i})=\{Q_{i1},Q_{i2}\}$, $i=1,...,5$ and we
set $\rho_{2}^{-1}(R_{s})=\{R_{s1},R_{s2}\}$,
$\rho_{2}^{-1}(A_{r})=\{A_{r1},A_{r2}\}$  where $s=1,...,2m$, $r=1,2$. 
Then 
$\langle i_{E}\rangle=G_{Q_{11}}=G_{Q_{12}}$,
$\langle i_{E^{'}}\rangle=G_{Q_{51}}=G_{Q_{52}}$,
$\langle i_{E}i_{E^{'}}\rangle=G_{Q_{i1}}=G_{Q_{i2}}$, $i=2,3,4$. Now it easy to
see that in $I)$ $\langle i_{{\widetilde{E}}^{'}}\rangle=G_{A_{rt}}$,
$r,t=1,2$, $\langle i_{{\widetilde{E}}}\rangle=G_{R_{s1}}=G_{R_{s2}}$, $s=1,...,2m$
which implies $t=8+8m$. In
$II)$:
$\langle i_{{\widetilde{E}}^{'}}\rangle=G_{A_{2t}}$, $t=1,2$,
$\langle i_{{\widetilde{E}}}\rangle=G_{R_{s1}}=G_{R_{s2}}$, $s=2,...,2m$,
$\langle i_{{\widetilde{E}}^{'}}i_{{\widetilde{E}}}\rangle=G_{R_{11}}=G_{R_{12}}$,
$t=12+8m$. In $III)$: 
$\langle i_{{\widetilde{E}}}\rangle=G_{R_{s1}}=G_{R_{s2}}$, $s=3,...,2m$, $\langle
i_{{\widetilde{E}}^{'}}i_{{\widetilde{E}}}\rangle
=G_{R_{11}}=G_{R_{12}}=G_{R_{21}}=G_{R_{22}}$,
$t=16+8m$.\qed


\section{$g(F)=3$}

The classification of infinite families of $G$-sandwich
canonically $3$-fibred surfaces is more difficult than the  genus -$2$ case. Moreover
there are many new cases with a rich structure. The following lemma is the analogue
of \ref{tredici}:

\begin{lemma}\label{diciannove}
Let $G$ be an Abelian group acting on a smooth genus-3 curve $F$ with quotient
$A$. Let $\eta\in H^{0}(F,\Omega^{1}_{F})$ be a
$1$-form such that ${\rm{div}}(\eta)$ is a $G$-invariant divisor and let
$H^{0}(F,\Omega^{1}_{F})=
\oplus_{\chi\in G^{\star}}V_{1,\chi}$ be the
decompositions into subspaces where $G$ acts on $V_{1,\chi}$ via the character
$\chi$. Then the following cases occur: 
\begin{equation}\label{tabellauno}
\begin{tabular}{|l|c|l|}\hline
  $g(A)$ &  $G$  & $H^{0}(F,\Omega^{1}_{F})$\\ \hline
$2$ & ${\bbbz}/2$ & $V_{1,{\rm{id}}}\oplus V_{1,\chi}$\\ \hline

$1$ & ${\bbbz}/2$ & $V_{1,{\rm{id}}}\oplus V_{1,\chi}$\\ \hline

$1$ & ${\bbbz}/3$ & $V_{1,{\rm{id}}}\oplus V_{1,\chi}\oplus V_{1,\chi^{2}}$\\ \hline

$1$ & ${\bbbz}/4$ & $V_{1,{\rm{id}}}\oplus V_{1,\chi}\oplus V_{1,\chi^{3}}$\\ \hline

$1$ & ${\bbbz}/2\times {\bbbz}/2$ & $V_{1,{\rm{id}}}\oplus V_{1,\chi_{1}}\oplus
V_{1,\chi_{2}}$\\ \hline

$0$ & ${\bbbz}/2$ & $V_{1,\chi}$ \\ \hline

$0$ & ${\bbbz}/3$ & $V_{1,\chi}\oplus V_{1,\chi^{2}}$\\ \hline

$0$ & ${\bbbz}/4$ & $V_{1,\chi}\oplus V_{1,\chi^{2}}\oplus V_{1,\chi^{3}}$\\ \hline

$0$ & ${\bbbz}/4$ & $V_{1,\chi}\oplus V_{1,\chi^{3}}$\\ \hline

$0$ & ${\bbbz}/2\times {\bbbz}/2$ & $V_{1,\chi_{1}}\oplus
V_{1,\chi_{1}\chi_{2}}\oplus V_{1,\chi_{2}}$\\ \hline 

$0$ & ${\bbbz}/2\times {\bbbz}/2$ & $V_{1,\chi_{1}}\oplus V_{1,\chi_{2}}$\\
\hline 

$0$ & ${\bbbz}/6$ & $V_{1,\chi}\oplus V_{1,\chi^{3}}\oplus V_{1,\chi^{5}}$\\ \hline

$0$ & ${\bbbz}/6$ & $V_{1,\chi}\oplus V_{1,\chi^{4}}\oplus V_{1,\chi^{5}}$\\ \hline

$0$ & ${\bbbz}/7$ & $V_{1,\chi}\oplus V_{1,\chi^{2}}\oplus V_{1,\chi^{4}}$\\ \hline

$0$ & ${\bbbz}/7$ & $V_{1,\chi}\oplus V_{1,\chi^{2}}\oplus V_{1,\chi^{5}}$\\ \hline

$0$ & ${\bbbz}/8$ & $V_{1,\chi}\oplus V_{1,\chi^{2}}\oplus V_{1,\chi^{5}}$\\ \hline

$0$ & ${\bbbz}/8$ & $V_{1,\chi}\oplus V_{1,\chi^{2}}\oplus V_{1,\chi^{3}}$\\ \hline

$0$ & ${\bbbz}/2\times {\bbbz}/4$ & $V_{1,\chi_{1}\chi_{2}}\oplus
V_{1,\chi_{1}\chi_{2}^{2}}\oplus V_{1,\chi_{2}}$\\ \hline 

$0$ & ${\bbbz}/2\times {\bbbz}/4$ & $V_{1,\chi_{1}\chi_{2}}\oplus
V_{1,\chi_{1}\chi_{2}^{2}}\oplus V_{1,\chi_{2}^{3}}$\\ \hline 

$0$ & ${\bbbz}/2\times {\bbbz}/2\times {\bbbz}/2$ & $V_{1,\chi_{100}}
\oplus V_{1,\chi_{010}}\oplus V_{1,\chi_{001}}$\\
\hline

$0$ & ${\bbbz}/9$ & $V_{1,\chi}\oplus V_{1,\chi^{2}}\oplus V_{1,\chi^{4}}$\\ \hline
 
$0$ & ${\bbbz}/12$ & $V_{1,\chi}\oplus V_{1,\chi^{2}}\oplus V_{1,\chi^{5}}$\\ \hline
 
$0$ & ${\bbbz}/14$ & $V_{1,\chi}\oplus V_{1,\chi^{3}}\oplus
V_{1,\chi^{5}}$\\ \hline

$0$ & ${\bbbz}/2\times {\bbbz}/8$ &
$V_{1,\chi_{2}^{3}}\oplus V_{1,\chi_{1}\chi_{2}^{2}}\oplus
V_{1,\chi_{2}}$\\ \hline

$0$ & ${\bbbz}/4\times {\bbbz}/4$ & $V_{1,\chi_{1}}\oplus
V_{1,\chi_{1}^{2}\chi_{2}^{3}}\oplus V_{1,\chi_{2}}$\\ \hline
 
\end{tabular}\end{equation}
\noindent
where $\chi$ is a generator of $({\bbbz}/d)^{\star}$ if
$d=2,3,4,6,7,8,9,12$, $\chi_{1}$,
$\chi_{2}$ generate 
$({\bbbz}/a\times {\bbbz}/b)^{\star}$ if $(a,b)=(2,2),(2,4),(2,7),(2,8),(4,4)$ and
$\chi_{100}$, $\chi_{010}$, $\chi_{001}$ generate $({\bbbz}/2\times {\bbbz}/2\times
{\bbbz}/2)^{\star}$.
\end{lemma}
\proof
The proof is straightforward but long. It requires two basic ingredients: a
standard Hurwitz's argument, plus the realization of $G$ as a $\Gamma$-quotient via a
homomorphism with torsion free kernel, where $\Gamma=\langle
a_{1},b_{1},...,a_{b},b_{b},x_{1},...,x_{m}\mid
\prod_{i=1}^{b}[a_{i},b_{i}]x_{1} \cdots x_{m}=x_{1}^{e_{1}}=\cdots
x_{m}^{e_{m}}=1\rangle$, $0\leq b\leq 2$, ${e_{i}}$, is the ramification index of each
$P_{ij}\in\pi^{-1}(Q_{i})$ and $Q_{1},...,Q_{m}\in A$ are the branch points
of $\pi:F\rightarrow A=F/G$.\qed\\
\\
We consider a surface $S$ with $p_{g}\gg 0$ and genus-$3$ fibre of the canonical
pencil $f:S\rightarrow B$. We call $b$ the genus of
$B$ and $q=q(S)$ the irregularity of $S$. From \ref{undici} we have to consider the
last four cases in the table (\ref{tabellacasi}): $A^{(3)}_{0}$ i.e. $q=b=0$;
$A^{(3)}_{1}$ i.e.
$q=1$, $b=0$; $A^{(3)}_{2}$ i.e. $q=2$, $b=0$; $B^{(3)}$ i.e. $q=b=1$. We divide our
classification-proof into two parts: $\mid G\mid \leq 8$ and $\mid G\mid >8$

\subsection{$\mid G\mid \leq 8$}
The examples of surfaces with isotrivial canonical
pencil with genus-$3$ fibre and $p_{g}\gg 0$ existing in the literature, as fas as we
know, are $G$-sandwich canonically $3$-fibred surfaces with a 
${\bbbz}/2$-action or a ${\bbbz}/2\times{\bbbz}/2$-action. As a by-product of our
classification we will see that all  $G$-sandwich canonically $3$-fibred surfaces
with 
$p_{g}\gg 0$ and $\mid G\mid \leq 8$ have $G={\bbbz}/2$ or
$G={\bbbz}/2\times{\bbbz}/2$ or $G={\bbbz}/2\times{\bbbz}/2\times{\bbbz}/2$.


\paragraph{\mbox{\bf{q=b=0, g=3}}:}
The surfaces constructed in this section are all new. In fact the case $q=b=0$ is the
richer and the harder one.
\begin{proposition}\label{qugualebugualezero}
If $S$ is a $G$-sandwich canonically $3$-fibred surface on a rational curve
with $\mid G\mid \leq 8$, $q(S)=0$, $p_{g}(S)=m-1\geq 3$ then $S$ is the
minimal desingularization of
$X=Z/G$ where $G$ acts diagonally on $Z=F\times
D$, $g(F)=3$, $X$ has $t$ nodes and $S$ is in one of the following
classes:

\noindent
{\bf{Case}} $G={\bbbz}/2\times{\bbbz}/2$:
$$
\begin{array}{rcll} 
I) & g(D)=2m-1& K^{2}_{S}=8\chi(\cO_{S})&t=16\\
II) & g(D)=2m-2 & K^{2}_{S}=8\chi(\cO_{S})-12&t=24\\
III) & g(D)=2m-3& K^{2}_{S}=8\chi(\cO_{S})-16&t=32\\
\end{array}
$$
\noindent
{\bf{Case}} $G={\bbbz}/2\times{\bbbz}/2\times{\bbbz}/2$:
$$
\begin{array}{rcll} 
i) & g(D)=4m+1& K^{2}_{S}=8\chi(\cO_{S})& t=0\\
ii) & g(D)=4m-1 & K^{2}_{S}=8\chi(\cO_{S})-4& t=16\\
iii) & g(D)=4m-3& K^{2}_{S}=8\chi(\cO_{S})-8&t=32\\
iv) & g(D)=4m-5& K^{2}_{S}=8\chi(\cO_{S})-12&t=48\\
v) & g(D)=4m-7& K^{2}_{S}=8\chi(\cO_{S})-16&t=64\\
\end{array}
$$
\end{proposition}
\proof
We have to consider groups $G$ in table (\ref{tabellauno}) with
$\mid G\mid\leq 8$ and $g(A)=0$. In fact we have to classify all the $G$
coverings
$D\rightarrow D/G=B\simeq I\!\!P^{1}$ such
that the induced
decomposition $H^{0}(D,\Omega^{1}_{D})=\oplus_{\chi\in
G^{\star}}V_{2,\chi}$ satisfies the two conditions in \ref{otto}. Now $G={\bbbz}/2$
does not occur since ${\rm{dim}}V_{1,{\rm{id}}}=0$, ${\rm{dim}}V_{2,{\rm{id}}}=0$ and
${\rm{dim}}V_{1,{\rm{-id}}}=3$, ${\rm{dim}}V_{2,{\rm{-id}}}>1$. $G={\bbbz}/3$ does
not occur. In fact we have to solve the system (\ref{zetatre}) which becomes a
simple system with integer coefficients; we can assume
${\rm{dim}}V_{2,\chi^{2}}=0$ where $G^{\star}=\langle\chi\rangle$. Then it leads to
the following solution: $1={\rm{dim}}V_{2,\chi^{}}=g(D)=p_{g}(S)$;
$0={\rm{dim}}V_{2,\chi^{}}=g(D)=p_{g}(S)$; which corresponds to a surface not
of general type. By the general theory
of Abelian covers $D\rightarrow B=D/{\bbbz}/4$ is given through
${\cal{L}}_{\chi},{\cal{L}}_{\chi^{2}}, {\cal{L}}_{\chi^{3}}\in {\rm{Pic}}(B)$
and three effective divisors $D_{1}$, $D_{2}$, $D_{3}$ such that:

\begin{equation}
\label{ilcasozetaquattro}
\left\{ \begin{array}{rl}
4{\cal{L}}_{\chi}= &D_{1}+2D_{2}+3D_{3}\\
{\cal{L}}_{\chi^{2}} = & 2{\cal{L}}_{\chi}-D_{2}-D_{3}\\
{\cal{L}}_{\chi^{3}}=& 3{\cal{L}}_{\chi}-D_{2}-2D_{3}.\\
\end{array}\right.
\end{equation}
We have to consider the two ${\bbbz}/4$-action in the table (\ref{tabellauno}).
If $H^{0}(F,\Omega^{1}_{F})=V_{1,\chi}\oplus V_{1,\chi^{2}}\oplus
V_{1,\chi^{3}}$, then we can assume ${\rm{dim}}V_{2,\chi^{3}}=0$ that is
${\rm{deg}}({\cal{L}}_{\chi^{3}})=1$. Since $D\rightarrow D/G=B$
is connected we have ${\rm{deg}}({\cal{L}}_{\chi^{2}})>0$ and
${\rm{deg}}({\cal{L}}_{\chi^{1}})>0$. Then
$g(D)\leq 1$ i.e. $S$ is not of general type. If
$H^{0}(F,\Omega^{1}_{F})=V_{1,\chi}\oplus V_{1,\chi^{3}}$ with
${\rm{dim}}V_{1,\chi}=2$ then ${\rm{dim}}V_{2,\chi^{3}}=0$. Moreover 
${\rm{deg}}({\cal{L}}_{\chi^{2}})>0$. Then we
obtain $i)$ ${\rm{deg}}({\cal{L}}_{\chi^{2}})=1$,
${\rm{deg}}({\cal{L}}_{\chi^{}})=2$ or $ii)$ ${\rm{deg}}({\cal{L}}_{\chi^{2}})=2$,
${\rm{deg}}({\cal{L}}_{\chi^{}})=3$. Then $g(D)=3$ and the minimal
desingularization $S$ has $p_{g}=2$.\\\\
{\bf{The case}} $G={\bbbz}/2\times{\bbbz}/2$.\\
 We recall that the stabilizer $G(P)$ of a point $P\in Z$ is the trivial group
or $G(P)={\bbbz}/2$. In particular by \ref{barlow} the minimal resolution
$S\rightarrow X=Z/G$ is the blow up of the $t$ nodes of $X$ and $S$ is minimal. We
want three reduced divisors
$D_{1},D_{2},D_{3}$ and three line bundles on
$I\!\!P^{1}$, ${\cal{L}}_{\chi_{1}}$ ${\cal{L}}_{\chi_{2}}$
${\cal{L}}_{\chi_{1}\chi_{2}}$ which satisfy the system (\ref{zetadueperse}). Set
$m={\rm{deg}}({\cal{L}}_{\chi_{1}})$ and we easily see that we obtain the desired
solutions.\\ 
\\
If $G={\bbbz}/6={\bbbz}/3\times{\bbbz}/2=H\times K$,
$G^{\star}=\langle\chi\rangle$,
$H^{\star}=\langle\phi\rangle$, $K^{\star}=\langle\psi\rangle$ where
$\chi_{\mid H}=\phi$ and $\chi_{\mid K}=\psi$ by \cite{P}[Proposition 2.1 and 2.15]we
have to solve
\begin{equation}
\label{ilcasozetadueperzetatre}
\left\{ \begin{array}{l}
6{\cal{L}}_{\chi}= D_{\chi}+2D_{\phi}+4D_{\phi^{2}}+3D_{\psi}+5D_{\chi^{-1}}\\
{\cal{L}}_{\chi^{2}}=2{\cal{L}}_{\chi}-D_{\phi^{2}}-D_{\psi}-D_{\chi^{-1}}\\
{\cal{L}}_{\chi^{3}}=3{\cal{L}}_{\chi}-D_{\phi}-D_{\phi^{2}}-D_{\psi}-2D_{\chi^{-1}}\\
{\cal{L}}_{\chi^{4}}=4{\cal{L}}_{\chi}-D_{\phi}-2D_{\phi^{2}}-2D_{\psi}-3D_{\chi^{-1}}\\
{\cal{L}}_{\chi^{5}}=5{\cal{L}}_{\chi}-D_{\phi}-3D_{\phi^{2}}-2D_{\psi}-4D_{\chi^{-1}}\\
\end{array}\right.
\end{equation}
Looking to the two actions in \ref{tabellauno} we easily see that we can assume that 
$V_{2,\chi^{5}}=0$ that is ${\rm{deg}}({\cal{L}}_{\chi})=1$. Then it is obvious
that $p_{g}(S)$ is small. A brute computation shows that if 
$H^{0}(F,\Omega^{1}_{F})=V_{1,\chi}\oplus V_{1,\chi^{3}}\oplus
V_{1,\chi^{5}}$ then $H^{0}(D,\Omega^{1}_{D})=V_{2,\chi}\oplus
V_{2,\chi^{2}}\oplus V_{1,\chi^{4}}$ with ${\rm{dim}}V_{2,\chi}=2$
${\rm{dim}} V_{2,\chi^{2}}={\rm{dim}} V_{2,\chi^{4}}=1$,
$g(D)=4$, $p_{g}(S)=2$ and
if $H^{0}(F,\Omega^{1}_{F})=V_{1,\chi}\oplus V_{1,\chi^{2}}\oplus
V_{1,\chi^{5}}$ then $H^{0}(D,\Omega^{1}_{D})=V_{2,\chi}\oplus
V_{2,\chi^{2}}\oplus V_{1,\chi^{3}}$ with ${\rm{dim}}V_{2,\chi}={\rm{dim}}
V_{2,\chi^{2}}={\rm{dim}} V_{2,\chi^{3}}=1$,
$g(D)=3$, $p_{g}(S)=1$
or $H^{0}(D,\Omega^{1}_{D})=V_{2,\chi}\oplus V_{2,\chi^{2}}\oplus
V_{1,\chi^{3}}$ with ${\rm{dim}}V_{2,\chi}=2$
${\rm{dim}} V_{2,\chi^{2}}={\rm{dim}} V_{2,\chi^{3}}=1$,
$g(D)=4$, $p_{g}(S)=2$.\\
\\
The cases $G={\bbbz}/7$ or $G={\bbbz}/8$ or $G={\bbbz}/2\times{\bbbz}/4$ have no
solutions.
\\
\\
{\bf{The case}} $G={\bbbz}/2\times{\bbbz}/2\times{\bbbz}/2$:\\
Let $G=H_{100}\times H_{010}\times H_{001}$,
$G^{\star}=\{\chi_{ijk}\}\cup\{{\rm{id}}\}$ where $H_{ijk}\simeq{\bbbz}/2$ are the
$7$ non-trivial subgroups and $\chi_{ijk}$ are the non-trivial characters of $G$ with
the obvious notation. We want divisors
$D_{ijk}$ and line bundles $L_{ijk}$ on $I\!\!P^{1}$ such that

\begin{equation}
\label{ilcasozetadueperzetadueperzetadue}
\left\{ \begin{array}{rl}
2{\cal{L}}_{100}= &D_{100}+D_{110}+D_{101}+D_{111}\\
2{\cal{L}}_{010}=&D_{010}+D_{110}+D_{011}+D_{111}\\
2{\cal{L}}_{001}=&D_{001}+D_{101}+D_{011}+D_{111}\\
{\cal{L}}_{110}=&{\cal{L}}_{100}+{\cal{L}}_{010}-D_{110}-D_{111}\\
{\cal{L}}_{101}=&{\cal{L}}_{100}+{\cal{L}}_{001}-D_{101}-D_{111}\\
{\cal{L}}_{011}=&{\cal{L}}_{010}+{\cal{L}}_{001}-D_{011}-D_{111}\\
{\cal{L}}_{111}=&{\cal{L}}_{100}+{\cal{L}}_{010}+{\cal{L}}_{001}-D_{110}-
D_{101}-D_{011}-D_{111}\\
\end{array}\right.
\end{equation}

Set $d^{1}_{ijk}={\rm{deg}}D_{ijk}$ in the case of a $G$-action on $F$ and 
$d^{2}_{ijk}={\rm{deg}}D_{ijk}$ in the case of a $G$-action on $D$. We recall
that $d^{s}_{ijk}$ is the number of points $Q$ on $C_{s}$ such that
the stabilizer $G_{P}=H_{ijk}$ where $P\in\pi_{s}^{-1}(Q)$ and
$\pi_{s}:C_{s}^{'}\rightarrow C_{s}$ is the
${\bbbz}/2\times{\bbbz}/2\times{\bbbz}/2$-Galois cover, $s=1,2$. The
$G$-action on a genus-3 curve $F$ with
$\pi:F\rightarrow F/G=A=I\!\!P^{1}$ is unique and it is described
by the following numbers:
$d^{1}_{100}=d^{1}_{010}=d^{1}_{001}=0$ $d^{1}_{110}=d^{1}_{011}=d^{1}_{101}=1$,
$d^{1}_{111}=2$. By table (\ref{tabellauno}) and theorem \ref{otto} we can assume
${\rm{deg}}{\cal{L}}_{010}= {\rm{deg}}{\cal{L}}_{001}=1$ in the system
(\ref{ilcasozetadueperzetadueperzetadue}). Moreover since $D$ is connected
${\rm{deg}}{\cal{L}}_{ijk}\geq 1$ for every triple $(i,j,k)$. Then we can easily
solve (\ref{ilcasozetadueperzetadueperzetadue}) and the desired $G$-actions on a curve
$D$ are:
$$
\begin{array}{rlll}\label{tabellanodi}
i)&d^{2}_{111}=d^{2}_{011}=d^{2}_{101}=d^{2}_{110}=0,&
d^{2}_{010}=d^{2}_{001}=2,& d^{2}_{100}=2m\\

ii)&d^{2}_{111}=d^{2}_{011}=d^{2}_{101}=0,&
d^{2}_{110}=d^{2}_{010}=1, d^{2}_{001}=2,& d^{2}_{100}=2m-1\\

iii)&d^{2}_{111}=d^{2}_{011}=0,& d^{2}_{101}=d^{2}_{110}=d^{2}_{010}=d^{2}_{001}=1,
&d^{2}_{100}=2m-2\\

iv)&d^{2}_{111}=d^{2}_{010}=d^{2}_{001}=0,&
d^{2}_{110}=d^{2}_{101}=d^{2}_{011}=1,&d^{2}_{100}=2m-2\\

v)&d^{2}_{111}=d^{2}_{010}=d^{2}_{001}=d^{2}_{011}=0,&
d^{2}_{110}=d^{2}_{101}=2,&d^{2}_{100}=2m-4.\dashv \\

\end{array}
$$                                       

\subsubsection{ A geometrical construction of the families of
\ref{qugualebugualezero} }\label{ioio} We will show a geometrical procedure to
construct the families found in
\ref{qugualebugualezero}.\\\\ {\bf{The case}}
$G={\bbbz}/2\times{\bbbz}/2$\label{rosso}:\\ we start with
$A=B=I\!\!P^{1}$. Let $\pi :C\rightarrow I\!\!P^{1}$ the
$2:1$-cover induced by the hyperelliptic involution $j_{1}:C\rightarrow C$ on a
genus-2 curve and let $W_{1},...,W_{6}\in C$ be the Weierstrass points. Set
$\pi(W_{i})=A_{i}$ $i=1,...,6$. Let $E\rightarrow I\!\!P^{1}$ the $2:1$-cover
branched on $A_{1}, A_{2}, A_{3}, A_{4}$ and $i_{1}:E\rightarrow E$ the induced
involution on the elliptic curve $E$.  The group $\langle i_{1}\rangle\times\langle
j_{1}\rangle\simeq{\bbbz}/2\times{\bbbz}/2$ acts on the
fibre product
$C\times_{I\!\!P^{1}}E=F$ in the obvious way, $\Gamma=F/\langle
i_{1}j_{1}\rangle$ is the $2:1$-cover
$\Gamma\rightarrow I\!\!P^{1}$ branched on $A_{5}, A_{6}$, $F/\langle
j_{1}\rangle=E$, $F/\langle
i_{1}\rangle=C$. Let $G^{\star}=\{ {\rm{id}},\chi_{1},\chi_{2},\chi_{3}\}$ where
$\chi_{3}=\chi_{1}\circ\chi_{2}$
$\chi_{1}(i_{1})=-1$, $\chi_{1}(j_{1})=1$, $\chi_{2}(i_{1})=1$,
$\chi_{2}(j_{1})=-1$ then
$H^{0}(F,\Omega^{1}_{F})=V_{1,\chi_{1}}\oplus V_{1,\chi_{2}}$ where  
$V_{1,\chi_{1}}=\langle \eta\rangle\simeq H^{0}(E,\Omega^{1}_{E})$, 
$V_{1,\chi_{2}}=\langle \mu_{1},\mu_{2}\rangle\simeq H^{0}(C,\Omega^{1}_{C})$. 
Let
${\widetilde{\pi}} :{\widetilde{C}}\rightarrow I\!\!P^{1}$,
${\widetilde{E}}\rightarrow I\!\!P^{1}$ be two $2:1$-cover branched respectively
on $R_{1},...,R_{2m}$ and on $Q_{1},Q_{2}$; then
$g({\widetilde{C}})=m-1$, $g({\widetilde{E}})=0$ . Let
$i_{2}:{\widetilde{C}}\rightarrow {\widetilde{C}}$,
$j_{2}:{\widetilde{E}}\rightarrow{\widetilde{E}}$ the induced involutions. The group  
$\langle i_{2}\rangle\times\langle j_{2}\rangle$ acts on the fibre product 
${\widetilde{C}}\times_{I\!\!P^{1}}{\widetilde{E}}=D$ and
${\widetilde{\Gamma}}=F/\langle i_{2}j_{2}\rangle$ is the $2:1$-cover
${\widetilde{\Gamma}}\rightarrow I\!\!P^{1}$ branched on $Q_{1}, Q_{2},
R_{1},...R_{2m}$, then $g({\widetilde{\Gamma}})=m$
and $D\rightarrow{\widetilde{\Gamma}}$ is unramified. In particular
$H^{0}(D,\Omega^{1}_{D})=V_{2,\chi_{1}}\oplus V_{2,\chi_{3}}$ where  
$V_{2,\chi_{1}}=\langle \alpha_{1},...\alpha_{m-1}\rangle\simeq
H^{0}({\widetilde{C}},\Omega^{1}_{{\widetilde{C}}})$, 
$V_{2,\chi_{3}}=\langle \beta_{1},...,\beta_{m}\rangle\simeq
H^{0}({\widetilde{\Gamma}},\Omega^{1}_{{\widetilde{\Gamma}}})$.
Consider $Z=F\times D$ with the $G$-action given by $\langle
i,j\rangle$ where $i=(i_{1},i_{2})$, $j=(j_{1},j_{2})$. 
In particular on $X=Z/G$ there are the $16$ nodes obtained by the $16$ points of $Z$
with stabilizer $\langle j\rangle$. Since
$H^{0}(X,\Omega^{2}_{X})=H^{0}(Z,\Omega^{2}_{Z})^{G}=[(V_{1,\chi_{1}}\oplus
V_{1,\chi_{2}})\otimes V_{2,\chi_{1}}\oplus V_{2,\chi_{3}})]^{G}$ then
$H^{0}(X,\Omega^{2}_{X})=V_{1,\chi_{1}}\otimes
V_{2,\chi_{1}})=\{\eta\wedge\alpha_{1},...,\eta\wedge\alpha_{m-1}\}$; that is $X$ is
the family $I)$ of \ref{qugualebugualezero}. The families $II)$ and $III)$
are obtained by a degeneration argument. In fact if the point $R_{1}$ coincides with
$Q_{1}$ we have $II)$ while $III)$ is given by the conditions: $R_{1}=Q_{1}$ and
$R_{2}=Q_{2}$.\\\\
{\bf{The case}} $G={\bbbz}/2\times{\bbbz}/2\times{\bbbz}/2$\\
We consider the ${\bbbz}/2\times{\bbbz}/2$-action on a genus-$2$ curve $C$
which gives the three families of  $G$-sandwich canonically $2$-fibred surfaces on a
rational curve with $q(S)=0$ described in section \ref{stanco}. With the same
notation there let
$\rho_{E^{''}}:E^{''}\rightarrow A=I\!\!P^{1}$ be the ${\bbbz}/2$-cover branched
on $Q_{2}$, $Q_{3}$ and let $i_{E^{''}}:E^{''}\rightarrow E^{''}$ be the corresponding
involution. The group $\langle i_{E^{''}}\rangle\times\langle
i_{E}\rangle\times\langle i_{E^{'}}\rangle$ acts on the normalization of the fibre
product $F=E^{''}\times_{A}C=E^{''}\times_{A}E\times_{A}E^{'}$
and let $\pi_{F}F\rightarrow A$ be
the ${\bbbz}/2\times{\bbbz}/2\times{\bbbz}/2$-Galois cover. We
set $\pi_{F}^{-1}(Q_{i})=\{Q_{i1}, Q_{i2}, Q_{i3},Q_{i4}\}$, in
particular $\pi_{F}^{\star}(Q_{i})=2Q_{i1}+2Q_{i2}+2Q_{i3}+2Q_{i4}$. It is
easy to see that $C=F/\langle i_{E^{''}}\rangle$, $C^{'}=F/\langle
i_{E}\rangle=E^{''}\times_{A}E$ and $C^{''}=F/\langle
i_{E^{'}}\rangle=E^{''}\times_{A}E^{'}$ are genus-$2$ curves and that
$Q_{1}=D^{1}_{101}$, $Q_{2}+Q_{3}=D^{1}_{111}$, $Q_{4}=D^{1}_{011}$,
$Q_{5}=D^{1}_{110}$. Moreover we can find a basis $\{\eta,\eta^{'},\eta^{''}\}$
of $H^{0}(F,\Omega^{1}_{F})$ such
that
$\langle\eta^{'},\eta^{''}\rangle=H^{0}(C,\Omega^{1}_{C})$,
$\langle\eta,\eta^{''}\rangle=H^{0}(C^{'},\Omega^{1}_{C^{'}})$,
$\langle\eta,\eta^{'}\rangle=H^{0}(C^{''},\Omega^{1}_{C^{''}})$.
Then $V_{1,\chi_{100}}=\langle\eta\rangle$,
$V_{1,\chi_{010}}=\langle\eta^{'}\rangle$,
$V_{1,\chi_{001}}=\langle\eta^{''}\rangle$. We have described
the ${\bbbz}/2\times{\bbbz}/2\times{\bbbz}/2$-action on $F$. We will construct
the suitable $G$-action on $D$. Let $B\rightarrow B$, $B^{'}\rightarrow
B$, $B^{''}\rightarrow B$ be the 
${\bbbz}/2$-covers of $B=I\!\!P^{1}$ branched on
$\{A^{010}_{1},A^{010}_{2}\}$, $\{A^{001}_{1},A^{001}_{2}\}$
and $\{A^{100}_{1},...,A^{100}_{2m}\}$ respectively and $i_{100}:B^{''}\rightarrow
B^{''}$, $i_{010}:B^{}\rightarrow B^{}$, $i_{001}:B^{'}\rightarrow B^{'}$ the
corresponding involutions. Let $D=B^{''}\times_{B}B\times_{B}B^{'}$
then
$G=\langle i_{100}, i_{010},i_{001}\rangle$ acts on $D$. We set 
$C_{100}=D/\langle i_{100}\rangle$, $C_{010}=D/\langle
i_{010}\rangle$, $C_{001}=D/\langle i_{001}\rangle$. Then 
$C_{100}$ is elliptic and $C_{100}\rightarrow B$ is branched on
$\{A^{010}_{1},A^{010}_{2}, A^{001}_{1},A^{001}_{2}\}$ while
$g(C_{010})=g(C_{001})=2m-1$ and $C_{010}\rightarrow B$, $C_{001}\rightarrow
B$ have branch locus
respectively $\{A^{001}_{1},A^{001}_{2}\}\cup\{A^{100}_{1},...,A^{100}_{2m}\}$
and $\{A^{010}_{1},A^{010}_{2}\}\cup\{A^{100}_{1},...,A^{100}_{2m}\}$. We can find a
basis 
$\langle\epsilon\rangle\cup\langle\alpha_{i}\rangle_{i=1,...,m-1}\cup
\langle\beta_{j}\rangle_{j=1,...,m}
\cup\langle\gamma_{r}\rangle_{r=1,...,m}\cup
\langle\delta_{s}\rangle_{s=1,...,m+1}$ of $H^{0}(D,\Omega^{1}_{D})$
such that $G$ acts in the following way:
\begin{equation}\label{azionecomplessa}
\begin{tabular}{|l|c|l|l|}\hline
  $$ &  $i_{100}^{\star}$  & $i_{010}^{\star}$& $i_{001}^{\star}$\\ \hline
$\epsilon$ & $\epsilon$ & $-\epsilon$& $-\epsilon$\\ \hline

$\alpha$ & $-\alpha$ & $\alpha$&$\alpha$\\ \hline

$\beta$ & $-\beta$ & $\beta$ & $-\beta$\\ \hline

$\gamma$ & $-\gamma$ &$-\gamma$ & $\gamma$\\ \hline

$\delta$ & $-\delta$ &$-\delta$ &$-\delta$\\ \hline
\end{tabular}
\end{equation}
\noindent in particular $V_{2,\chi_{100}}=\langle\alpha_{1},...,\alpha_{m-1}\rangle$,
$V_{2,\chi_{010}}=0$, $V_{2,\chi_{001}}=0$. Now it is easy to construct a
$G={\bbbz}/2\times{\bbbz}/2\times{\bbbz}/2$-action on the product $Z=F\times
D$ such that $X=Z/G$,
$H^{0}(X,\Omega^{2}_{X})=H^{0}(Z,\Omega^{2}_{Z})^{G}=V_{1,\chi_{100}}\otimes
V_{2,\chi_{100}})=\{\eta\wedge\alpha_{1},...,\eta\wedge\alpha_{m-1}\}$ and $S=X/G$ is
the surface with $p_{g}(S)=m-1$ of type $i)$ in \ref{qugualebugualezero}. We can
obtain the other surfaces in \ref{qugualebugualezero} with a
${\bbbz}/2\times{\bbbz}/2\times{\bbbz}/2$-action through a degeneration
argument similar to the previous one in the ${\bbbz}/2\times{\bbbz}/2$-case.

\paragraph{\mbox{\bf{q=1, b=0, g=3}}:} as far as we know the only examples of
isotrivial canonical fibrations with $q=1$, $b=0$, $g=3$ are in \cite{Z2}, see also
\cite{Ca} and \cite{Ca1}. Here we will give the complete classification of
$G$-sandwich canonically
$3$-fibred surfaces on a rational curve with $p_{g}(S)\gg 0$, $q=1.$ 

\begin{proposition}\label{pippo}
If $S$ is a $G$-sandwich canonically $3$-fibred surface on a rational curve
with $q(S)=1$, $p_{g}(S)=m-1\geq 2$ then $S$ is the minimal desingularization of
$X=Z/G$ where $G={\bbbz}/2\times{\bbbz}/2$ acts diagonally on $Z=F\times
D$, $g(F)=3$ and $S$ is in one of the following classes:
$$
\begin{array}{rcll} 
I) & g(D)=2m-1& K^{2}_{S}=8\chi(\cO_{S})&t=0\\
II) & g(D)=2m-2 & K^{2}_{S}=8\chi(\cO_{S})-4&t=8\\
III) & g(D)=2m-3& K^{2}_{S}=8\chi(\cO_{S})-8&t=16\\
\end{array}
$$
\noindent
\end{proposition}
\proof It is easier than the proof of \ref{qugualebugualezero}. In fact
$g(A)=1$ and by the table (\ref{tabellauno}) we have to classify all the
$G$-coverings $D\rightarrow D/G=B\simeq I\!\!P^{1}$ where
$G={\bbbz}/2$ or 
$G={\bbbz}/3$ or $G={\bbbz}/4$ or $G={\bbbz}/2\times{\bbbz}/2$ such that the
induced decomposition $H^{0}(D,\Omega^{1}_{D})=\oplus_{\chi\in
G^{\star}}V_{2,\chi}$ satisfies the two conditions in \ref{otto}. 
We exclude $G={\bbbz}/2$, $G={\bbbz}/3$, $G={\bbbz}/4$ by a simple computation as
in \ref{qugualebugualezero}. Assume $G={\bbbz}/2\times{\bbbz}/2$. We want three
reduced divisors $D_{1},D_{2},D_{3}$ and three line bundles on
$I\!\!P^{1}$, ${\cal{L}}_{\chi_{1}}$, ${\cal{L}}_{\chi_{2}}$,
${\cal{L}}_{\chi_{1}\chi_{2}}$ under the
constraints
${\rm{deg}}({\cal{L}}_{\chi_{2}})=1$,
${\rm{deg}}({\cal{L}}_{\chi_{1}\chi_{2}})\geq1$. Set
$m={\rm{deg}}({\cal{L}}_{\chi_{1}})$ and we easily solve the
system (\ref{zetadueperse}) which gives the desired solutions\qed

\subsubsection{The three families of proposition \ref{pippo}: a geometrical
construction}
We want to construct a ${\bbbz}/2\times{\bbbz}/2$-cover $\pi_{1}:F\rightarrow
A$ on the elliptic curve $A$. Let $C\rightarrow A$ be the $2$-to-$1$
cover branched on $P,Q\in A$, $\pi_{1}^{-1}(P)=P_{1}$,
$\pi_{1}^{-1}(Q)=Q_{1}$ and set $C^{'}$ another copy of $C$; call $i_{C}:C\rightarrow
C$ and  $i_{C^{'}}:C^{'}\rightarrow C^{'}$ the associated involutions, let
$\nu:F\rightarrow C\times_{A}C^{'}$ be the normalization map and
set $\{P_{1}^{1},P_{1}^{2}\}=\nu^{-1}((P_{1},P_{1}))$,
$\{Q_{1}^{1},Q_{1}^{2}\}=\nu^{-1}((Q_{1},Q_{1}))$. Obviously
${\bbbz}/2\times{\bbbz}/2=\{i_{C},i_{C^{'}}\}$ acts on $F$ and it interchanges
$P_{1}^{1}$ with $P_{1}^{2}$ and $Q_{1}^{1}$ with $Q_{1}^{2}$. In
particular $\Gamma=F/\langle i_{1}\dot j_{1}\rangle$ is an elliptic curve,
$H^{0}(F,\Omega^{1}_{F})=V_{1,{\rm{id}}}\oplus V_{1,\chi_{1}}\oplus
V_{1,\chi_{2}}$, $V_{1,{\rm{id}}}\oplus
V_{1,\chi_{1}}=H^{0}(C,\Omega^{1}_{C})$, $V_{1,{\rm{id}}}\oplus
V_{1,\chi_{2}}=H^{0}(C^{'},\Omega^{1}_{C^{'}})$,
$V_{1,{\rm{id}}}=H^{0}(\Gamma,\Omega^{1}_{\Gamma})$. Now consider
$\pi_{2}:{\widetilde{C}}\rightarrow I\!\!P^{1}=B$ a double cover branched on
$R_{1},...,R_{2m}$ points and let ${\widetilde{C^{'}}}\rightarrow B$ be the
double cover branched on $A_{1},A_{2}$. Then the normalization
${\widetilde{\nu}}: D\rightarrow
{\widetilde{C}}\times_{B}{\widetilde{C}^{'}}$ is a
$2$-to-$1$ covering of ${\widetilde{C}}$ and $g(D)=2m-1$. If 
$i_{\widetilde{C}}:{\widetilde{C}}\rightarrow
{\widetilde{C}}$ and  $i_{\widetilde{C^{'}}}:{\widetilde{C^{'}}}\rightarrow
{\widetilde{C^{'}}}$ are the associated involutions then
${\bbbz}/2\times{\bbbz}/2=\{i_{\widetilde{C}},i_{\widetilde{C^{'}}}\}$ acts on
$D$ and
${\widetilde{\Gamma}}=D/\langle
i_{\widetilde{C^{'}}}i_{\widetilde{C^{'}}}\rangle$ is branched on
$R_{1},...,R_{2m},A_{1},A_{2}$ and it has $g({\widetilde{\Gamma}})=m$. Then with the
standard notations we have $H^{0}(D,\Omega^{1}_{D})=
V_{2,\chi_{1}}\oplus V_{2,\chi_{3}}$ where
$V_{2,\chi_{1}}=H^{0}({\widetilde{C}},\Omega^{1}_{{\widetilde{C}}})$,
$V_{2,\chi_{3}}=H^{0}({\widetilde{\Gamma}},\Omega^{1}_{{\widetilde{\Gamma}}})$. Let
$Z=F\times D$, let $G\times Z\rightarrow Z$ be the diagonal
action where $G=\langle (i_{C},i_{\widetilde{C}})\rangle\times\langle
(i_{{C^{'}}},i_{\widetilde{C^{'}}})\rangle$ and
$X=Z/G$. By the proof of \ref{otto} we see that $H^{0}(X,\Omega^{2}_{X})\simeq
V_{1,\chi_{1}}\otimes V_{2,\chi_{1}}$ that is
$\Phi_{\mid K_{X}\mid}$ is composed with the rational pencil $f_{2}:X\rightarrow
B$ and $X$ is in the class $I)$ of \ref{pippo}. The classes $II)$ and $III)$ can
be obtained by a degeneration argument from the class $I)$. In fact if $R_{1}=A_{1}$
we have $II)$ and if $R_{1}=A_{1}$, $R_{2}=A_{2}$ we have $III)$.


\paragraph{\mbox{\bf{q=b=1, g=3}}:} In \cite{X2} Xiao constructed an infinite
family of 
${\bbbz}/2\times{\bbbz}/2$-sandwich canonically $3$-fibred surfaces on an elliptic
curve. In the following proposition we will classify all the $G$-sandwich canonically 
$3$-fibred surfaces with $q=b=1$, $p_{g}(S)\gg 0$. In particular
we will show that $a)$ the family of Xiao is the only
${\bbbz}/2\times{\bbbz}/2$-infinite family; $b)$ there is also a
${\bbbz}/2\times{\bbbz}/2\times\bbbz/2$-infinite family.

\begin{proposition}\label{pippodue}
If $S$ is a $G$-sandwich canonically $3$-fibred surface on an elliptic curve
with $p_{g}(S)=m\geq 2$ then $q(S)=1$, $S=Z/G$ where $G$ acts diagonally on
$Z=F\times D$, $g(F)=3$ and $S$ is in one of the
following classes:

\noindent
{\bf{Case}} $G={\bbbz}/2\times{\bbbz}/2$:
$$
\begin{array}{rcll} 
I) & g(D)=2m+1& K^{2}_{S}=8\chi(\cO_{S}). & \\
\end{array}
$$
\noindent
{\bf{Case}} $G={\bbbz}/2\times{\bbbz}/2\times{\bbbz}/2$:
$$
\begin{array}{rcll} 
II) & g(D)=4m+1& K^{2}_{S}=8\chi(\cO_{S})& \\
\end{array}
$$
\end{proposition}
\proof
Since $q(S)=1$ and the fibre of the Albanese map is connected, we have that
$f_{2}:D\rightarrow B=D/G$ is the Albanese map $g(B)=1$ and
$g(A)=0$. By (\ref{tabellauno}) we have to consider all the $G$ actions on
$D$ such that the induced decomposition $H^{0}(D,\Omega^{1}_{D})=\oplus_{\chi\in
G^{\star}}V_{2,\chi}$ satisfies the two conditions in \ref{otto}. We can exclude all
the cases with the exceptions of 
$G={\bbbz}/2\times{\bbbz}/2$, $G={\bbbz}/2\times{\bbbz}/2\times{\bbbz}/2$ as in
\ref{qugualebugualezero}.
\\\\
{\bf{The case}} $G={\bbbz}/2\times{\bbbz}/2$.\\
We consider the system (\ref{zetadueperse}) adapted to
our hypothesis. By the first 
${\bbbz}/2\times{\bbbz}/2$-action on $F$ we obtain that
${\cal{L}}_{\chi_{1}}$, ${\cal{L}}_{\chi_{2}}$ and ${\cal{L}}_{\chi_{1}\chi_{2}}$ are
non-trivial, torsion line bundles, then $g(D)=1$, i.e. $S$ is not of general
type. Let us consider the other ${\bbbz}/2\times{\bbbz}/2$-action on $F$.
>From (\ref{zetadueperse}) we obtain that
${\cal{L}}_{\chi_{2}}$ is a non-trivial torsion line bundle on the elliptic curve
$B$ and
${\cal{L}}_{\chi_{1}}+{\cal{L}}_{\chi_{2}}\equiv{\cal{L}}_{\chi_{1}\chi_{2}}$ where
${\rm{deg}}{\cal{L}}_{\chi_{1}}=m$. Then we have an unique solution and
$g(D)=2m+1$.\\\\\
{\bf{The case $G={\bbbz}/2\times{\bbbz}/2\times{\bbbz}/2$.}}\\
Since $H^{0}(F,\Omega^{1}_{F})=V_{1,\chi_{100}}\oplus
V_{1,\chi_{010}}\oplus V_{1,\chi_{001}}$, we can
require $V_{2,\chi_{010}}=V_{2,\chi_{001}}=0$ and
${\rm{dim}}V_{2,\chi_{100}}=p_{g}(S)=m>1$. We maintain the notations of the
system (\ref{ilcasozetadueperzetadueperzetadue}) and it follows that
$L_{010}$,$L_{001}$ and
$L_{011}$ are three distinct, non-trivial $2$-torsion bundles on the elliptic curve
$B$. Moreover $L_{110}=L_{100}+L_{010}$,  $L_{101}=L_{100}+L_{001}$, 
$L_{111}=L_{100}+L_{011}$. Then
${\rm{dim}}V_{2,\chi_{100}}={\rm{dim}}V_{2,\chi_{110}}=
{\rm{dim}}V_{2,\chi_{101}}={\rm{dim}}V_{2,\chi_{111}}=m$ and $g(D)=4m+1$.\qed

\subsubsection{Geometrical construction of the two families of
proposition \ref{pippodue}}
{\bf{The case $G={\bbbz}/2\times{\bbbz}/2:$}}\\
The action $G\times F\rightarrow F$ has been described above;
see \ref{rosso}. Let ${\widetilde{\pi}} :{\widetilde{C}}\rightarrow B$,
$g(B)=1$ be a $2:1$-cover branched on $R_{1},...,R_{2m}$, let
${\widetilde{E}}\rightarrow B$ a $2$-to-$1$ unramified covering,
$g({\widetilde{E}})=1$ and we denote by
$i_{2}:{\widetilde{C}}\rightarrow {\widetilde{C}}$,
$j_{2}:{\widetilde{E}}\rightarrow{\widetilde{E}}$ the induced involutions. On the
fibre product ${\widetilde{E}}\times_{B}{\widetilde{C}}=D$ acts the
group $\langle i_{2}\rangle\times\langle j_{2}\rangle\equiv{\bbbz}/2\times{\bbbz}/2$
and ${\widetilde{\Gamma}}=D/\langle i_{2}\dot j_{2}\rangle$ is an unramified
$2$-to-$1$ cover ${\widetilde{\Gamma}}\rightarrow B$. We can easily
find a basis $\langle
\alpha_{1},...\alpha_{m};\beta_{1},...\beta_{m};\epsilon\rangle=
H^{0}(D,\Omega^{1}_{D})$ such
that $H^{0}({\widetilde{C}},\Omega^{1}_{{\widetilde{C}}})\simeq\langle
\alpha_{1},...\alpha_{m};\epsilon\rangle$,
$H^{0}({\widetilde{E}},\Omega^{1}_{{\widetilde{E}}})\simeq\langle\epsilon\rangle\simeq
H^{0}({\widetilde{\Gamma}},\Omega^{1}_{{\widetilde{\Gamma}}})\simeq$ and
$\langle\beta_{1},...\beta_{m};\epsilon\rangle$ is also $i_{2}$ anti-invariant .
Moreover 
$V_{2,\chi_{1}}=\langle\alpha_{1},...,\alpha_{m}\rangle$,
$V_{2,\chi_{3}}=\langle\beta_{1},...,\beta_{m}\rangle$,
$V_{2,{\rm{id}}}=\langle\epsilon\rangle$. Consider $Z=F\times D$
with the $G$-action given by $\langle i,j\rangle$ where $i=(i_{1},i_{2})$,
$j=(j_{1},j_{2})$ and let $S=Z/G$. The action is free and since
$H^{0}(S,\Omega^{2}_{S})=H^{0}(Z,\Omega^{2}_{Z})^{G}=[(V_{1,\chi_{1}}\oplus
V_{1,\chi_{2}})\otimes V_{2,\chi_{1}}\oplus V_{2,\chi_{12}})]^{G}$ then
$H^{0}(S,\Omega^{2}_{X})=V_{1,\chi_{1}}\otimes
V_{2,\chi_{1}})=\{\eta\wedge\alpha_{1},...,\eta\wedge\alpha_{m}\}$; that is
$\Phi_{\mid K_{S}\mid}$ yields a pencil with image the elliptic normal curve of
degree $m$.\\\\
{\bf{The case $G={\bbbz}/2\times{\bbbz}/2\times{\bbbz}/2:$}}\\
We have constructed the suitable
${\bbbz}/2\times{\bbbz}/2\times{\bbbz}/2$-action on  the genus-$3$ curve
$F$ in \ref{ioio}. It is easy to produce the action on
$D$. In fact Let $B^{''}\rightarrow B$ the $2$-to-$1$ cover of
the elliptic curve $B$ branched on $2m$ points
$\{A^{100}_{i}\}_{i=1,...2m}$ and let $B\rightarrow B$,
$B^{''}\rightarrow B$ be two distinct 
${\bbbz}/2$-unramified covers of $B$. Let $i_{100}:B^{''}\rightarrow
B^{''}$, $i_{010}:B^{}\rightarrow B^{}$, $i_{001}:B^{'}\rightarrow B^{'}$
be the corresponding involutions. If
$D=B^{''}\times_{B}B\times_{B}B^{'}$ then
$G=\langle i_{100}, i_{010},i_{001}\rangle$ acts on $D$. We set 
$C_{100}=D/\langle i_{100}\rangle$, $C_{010}=D/\langle
i_{010}\rangle$, $C_{001}=D/\langle i_{001}\rangle$. Then 
$C_{100}$ is elliptic, while
$g(C_{010})=g(C_{001})=2m+1$ and $C_{010}\rightarrow B$,
$C_{001}\rightarrow B$ have branch locus $\{A^{100}_{i}\}_{i=1,...2m}$.
It is now easy to see that the
decomposition $\oplus_{\chi\in
G^{\star}}V_{2,\chi}$ and the diagonal action on the product surface
$Z=F\times D$ satisfy \ref{otto} and they produce the
$G$-family of \ref{pippodue}.

\paragraph{\mbox{\bf{q=2, b=0, g=3}}:}in \cite{Be}[Example 2], Beauville constructed
an infinite family of surfaces with $q=2$, canonical map composed with a pencil and
non-surjective Albanese map. In \cite{K}[Theorem 3.6] Konno showed that this family
is essentially unique: that is if $S$ has canonical map composed with a pencil,
$q=2$, $p_{g}\geq 8$ and non-surjective Albanese map then $S$ is the example of
Beauville. In the following proposition we give a simple proof of Konno's result in
the case of 
$G$-sandwich canonically fibred surfaces. We remark that we does not assume $g=3$. 

\begin{proposition}\label{qugualedue}
If $S$ is a $G$-sandwich canonically fibred surface with $q(S)=2$, $p_{g}(S)=m\geq 2$
then $S=Z/G$ where $G=\bbbz/2$, the group acts diagonally and freely on $Z=F\times
D$, $g(F)=3$, $g(D)=m$, $F/G=A$ has genus $2$, the
base $D/G=B$ of the canonical fibration is rational and
$K_{S}^{2}=8\chi(\cO_{S})$.
\end{proposition}
\proof Assume that $S$ is a minimal surface with $q(S)=2$ and with canonical
fibration $f:S\rightarrow B$ with fibre $F$ of genus $g>1$. Then $K_{S}\equiv
Z+f^{\star}(\alpha)$ where
$Z$, $\alpha$ are effective divisors on $S$ and $B$
respectively and $h^{0}(B,\cO_{B}(\alpha))=p_{g}(S)$. By \ref{novedieci} $B$ is
rational. By the well-known Miyaoka-Yau inequality we have 
$$9(p_{g}-1)\geq
K_{S}^{2}=K_{S}Z+K_{S}f^{\star}(\alpha)\geq (2g-2)(p_{g}-1),
$$
then $2\leq g\leq 5$. In particular if $S$ is a $G$-sandwich canonically fibred
surface with top $Z=F\times D$ and base $Z=A\times B$
then $B=B$, $3\leq g(F)\leq 5$, $g(A)=2$. It is very easy to classify
all the Abelian actions with group $G$ such that $F/G=A$:

$$\begin{tabular}{|l|c|l|}\hline
  $g(F)$ &  $G$  & $H^{0}(F,\Omega^{1}_{F})$\\ \hline
$3$ & ${\bbbz}/2$ & $V_{1,{\rm{id}}}\oplus V_{1,\chi}$\\ \hline

$4$ & ${\bbbz}/2$ & $V_{1,{\rm{id}}}\oplus V_{1,\chi}$\\ \hline

$4$ & ${\bbbz}/3$ & $V_{1,{\rm{id}}}\oplus V_{1,\chi}\oplus V_{1,\chi^{2}}$\\ \hline

$5$ & ${\bbbz}/2$ & $V_{1,{\rm{id}}}\oplus V_{1,\chi}$\\ \hline

$5$ & ${\bbbz}/4$ & $V_{1,{\rm{id}}}\oplus V_{1,\chi}\oplus V_{1,\chi^{2}}\oplus
V_{1,\chi^{3}}$\\ \hline

$5$ & ${\bbbz}/2\times{\bbbz}/2$ & $V_{1,{\rm{id}}}\oplus V_{1,\chi_{1}}\oplus
V_{1,\chi_{2}}\oplus V_{1,\chi_{12}}$\\ \hline
\end{tabular}
$$
\noindent
Now we must classify all the $G$ actions on a curve $D$ such that
$D/G=B$ and the induced
decomposition $H^{0}(D,\Omega^{1}_{D})=\oplus_{\chi\in
G^{\star}}V_{2,\chi}$ satisfies \ref{otto}. It is an easy computation to show that
only the case $g(F)=3$ works and it gives Beauville's family.\qed\\

We have proved the first theorem stated in the introduction. Now we will show the
second one.

\subsection{$\mid G\mid>8$}
The results showed in this section are new. We will prove that there is a
rich geometry among the $G$-sandwich canonically $3$-fibred surfaces with $\mid
G\mid>8$. More precisely we will show that only the case $G={\bbbz}/2\times{\bbbz}/8$
occurs but it gives many different cases. 

\begin{theorem}\label{eccolotto}
If $S$ is a $G$-sandwich canonically $3$-fibred surface on a rational curve
with $\mid G\mid>8$ and $p_{g}(S)\geq 3$ then $q(S)=0$, $S$ is the minimal
desingularization of $X=Z/G$ where
$G={\bbbz}/2\times{\bbbz}/8$ acts diagonally on
$Z=F\times D$, $g(F)=3$ and $S$ is in one of the following
classes:

$$
\begin{array}{rlll} 
i) & g(D)=8m-1-4(a+b)& p_{g}(S)=m-a& 0\leq a\leq 1,  0\leq b\leq 1 \\

ii) & g(D)=8m-4(a+b-1)& p_{g}(S)=m+1-(a+b)& 0\leq a\leq 2,  0\leq b\leq 1 \\

iii) & g(D)=8m+3-4(a+b)& p_{g}(S)=m-a& 0\leq a\leq 1,  0\leq b\leq 1 \\

iv) & g(D)=8m+8-4(a-1)& p_{g}(S)=m+1-a& 0\leq a\leq 2 \\

\end{array}
$$
\noindent
where $a,b,m\in\bbbz^{+}$, $m\geq 4$. Moreover the action of $G$ is described 
completely.
\end{theorem}
\proof
We prove the claim by a direct computation. We need to compute all the $G$ actions on
a curve $D$ such that $B=D/G$ is a rational curve and the
induced decomposition $H^{0}(D,\Omega^{1}_{D})=\oplus_{\chi\in
G^{\star}}V_{2,\chi}$ satisfies \ref{otto}.
\\
By the standard theory of Abelian
covers [see: \cite{P}] we can write the branch
locus $\Delta$ of $\rho_{D}:D\rightarrow B$ in the following
way: $\Delta=\sum_{H\in\cC}\sum_{\psi\in S_{H}}D_{\psi,H}$ where $\cC$ is the set of
cyclic subgroups of $G$, $S_{H}$ is the set of generators of $H^{\star}$ for every
$H\in\cC$ and $D_{\psi,H}$ is the reduced sum of all points $P$ of $D$ such that $H$
is the stabilizer of a(ny) point in $\rho_{D}^{-1}(P)$ and $H$ operates via $\psi$ on
the cotangent fibre on $P$.
\\
{\bf{ The case}} $G=\bbbz/9$:\\
let $H={\bbbz}/3$ be the cyclic
proper subgroup of $G$ and let $\langle \chi\rangle= G^{\star}$, $\langle
\phi\rangle=H^{\star}$. The Galois morphism $\rho_{D}:D\rightarrow
B$ is given by a line bundle
$L=L_{\chi}$ on $B$ such that:  
$$
9L\equiv
D_{\chi^{}}+5D_{\chi^{2}}+3D_{\phi^{}}+6D_{\phi^{2}}
+7D_{\chi^{4}}+2D_{\chi^{5}}+4D_{\chi^{7}}+8D_{\chi^{8}}.
$$
By table (\ref{tabellauno}) we have to consider the following three cases:
\begin{enumerate}
\item[1)] ${\rm{dim}}V_{2,\chi^{-1}}=p_{g}(S)$ and
$V_{2,\chi^{-2}}=V_{2,\chi^{-4}}=0$;
\item[2)] ${\rm{dim}}V_{2,\chi^{-2}}=p_{g}(S)$ and
$V_{2,\chi^{-1}}=V_{2,\chi^{-4}}=0$;
\item[3)] ${\rm{dim}}V_{2,\chi^{-4}}=p_{g}(S)$ and
$V_{2,\chi^{-1}}=V_{2,\chi^{-2}}=0$.
\end{enumerate}
On the other hand all these cases differ by an automorphism of $G$. Then we can
compute the first action only: By \cite{P} we have:
$$ \label{eq:zetanoveperse}
\left\{ \begin{array}{ll}

9L& \equiv
D_{\chi^{}}+5D_{\chi^{2}}+3D_{\phi^{}}+6D_{\phi^{2}}
+7D_{\chi^{4}}+D_{\chi^{}}+2D_{\chi^{5}}+4D_{\chi^{7}}+8D_{\chi^{8}}\\

{\cal{L}}_{\chi^{-2}}&=7L-
3D_{\chi^{2}}-2D_{\phi^{}}-4D_{\phi^{2}}-5D_{\chi^{4}}-D_{\chi^{5}}-3D_{\chi^{7}}-6

D_{\chi^{8}}\\ 
{\cal{L}}_{\chi_{-4}} &=5L-
2D_{\chi^{2}}-D_{\phi^{}}-3D_{\phi^{2}}-3D_{\chi^{4}}-D_{\chi^{5}}-2D_{\chi^{7}}-4
D_{\chi^{8}}.\\
\end{array}\right.
$$
\noindent
Let $l_{i}={\rm{deg}}L_{\chi_{i}}$, $i=1,\dots,8$, $n_{i}={\rm{deg}}D_{\chi^{i}}$ and
$n_{3}={\rm{deg}}D_{\phi^{}}$, $n_{6}={\rm{deg}}D_{\phi^{2}}$. Taking degrees, by
\ref{eq:zetanoveperse} we obtain:

$$ \label{eq:zetanovepersegradi}
\left\{ \begin{array}{ll}
9l_{1}& =
n_{1}+5n_{2}+3n_{3}+6n_{6}+7n_{4}+2n_{5}+4n_{7}+8n_{8}\\

l_{7}& =7l_{1}-
(3n_{2}+2n_{3}+4n_{6}+5n_{4}+n_{5}+3n_{7}+6n_{8})\\

l_{5}& =5l_{1}-
(2n_{2}+n_{3}+3n_{6}+3n_{4}+n_{5}+2n_{7}+4n_{8})\\
\end{array}\right.
$$
\noindent
Since $g(B)=0$ then $l_{7}=l_{5}=1$ and this implies $n_{1}=0, n_{2}=n_{7}=1$ or
$n_{2}=0, n_{1}=n_{8}=1$; that is $g(D)=0$: a contradiction.
\\
\\
{\bf{ The case}} $G=\bbbz/12$:\\
let $H_{2}={\bbbz}/2$, $H_{3}={\bbbz}/3$, $H_{4}={\bbbz}/4$, $H_{6}={\bbbz}/6$ be the
cyclic proper subgroups of $G$ and let $\langle \chi\rangle= G^{\star}$, $\langle
\mu_{i}\rangle= H_{i}^{\star}$. To construct $\rho_{D}:D\rightarrow
B$ we need a line bundle $L=L_{\chi}$ on $B$ such that:
$$
12L\equiv
D_{\chi^{}}+5D_{\chi^{5}}+7D_{\chi^{7}}+11D_{\chi^{11}}+6D_{\mu_{2}}
4D_{\mu_{3}}+8D_{\mu_{3}^{2}}+3D_{\mu_{4}}+9D_{\mu_{4}^{3}}+2D_{\mu_{6}}+
10D_{\mu_{6}^{5}}.
$$
By \ref{otto} and table (\ref{tabellauno}) we can assume ${\rm{deg}}L=l=1$, then
there are not infinite families. However we can easily find all these surfaces. Set
$l_{i}={\rm{deg}}L_{\chi^{i}}$. The case
$l_{5}=1$ and $l_{2}\geq 3$ does not occur; the case
$l_{5}\geq 3$ and $l_{2}=1$ has two solutions: $i)$ ${\rm{deg}}D_{\chi^{}}=7$,
${\rm{deg}}D_{\chi^{5}}=1$,
$D_{\chi^{7}}=D_{\chi^{11}}=D_{\mu_{2}}=D_{\mu_{3}}=D_{\mu_{3}^{2}}=D_{\mu_{4}}=
D_{\mu_{4}^{3}}=D_{\mu_{6}}=D_{\mu_{6}^{5}}=0$; $ii)$
${\rm{deg}}D_{\chi^{}}=6$,
${\rm{deg}}D_{\mu_{2}}=1$,
${\rm{deg}}D_{\chi^{5}}=D_{\chi^{7}}=D_{\chi^{11}}=D_{\mu_{3}}=D_{\mu_{3}^{2}}=D_{\mu_{4}}=
D_{\mu_{4}^{3}}=D_{\mu_{6}}=D_{\mu_{6}^{5}}=0$; both solutions have $p_{g}(S)=2$.
\\
\\
{\bf{ The case}} $G=\bbbz/14$:\\ 
Let $H={\bbbz}/2$, $K={\bbbz}/7$ be the
cyclic proper subgroups of $G$ and let $\langle \chi\rangle= G^{\star}$, 
$\langle\mu\rangle= H^{\star}$, $\langle\psi\rangle= K^{\star}$. To construct
$\rho_{D}$ we need a line bundle $L=L_{\chi}$ such that:
$$
14L\equiv
D_{\chi^{}}+5D_{\chi^{3}}+3D_{\chi^{5}}+11D_{\chi^{9}}+9D_{\chi^{11}}+13D_{\chi^{13}}+
7D_{\mu}
+\sum_{i=1}^{6}2iD_{\chi^{2i}}.
$$
By \ref{otto} and table (\ref{tabellauno}) we can assume ${\rm{deg}}L=l=1$, then
there are not infinite families. Moreover both cases
$l_{5}=1$ and $l_{3}\geq 3$ and
$l_{5}\geq 3$ and $l_{3}=1$ have solutions, but $p_{g}(S)\leq 2$.
\\
\\
{\bf{ The case}} $G=\bbbz/4\times\bbbz/4$:
\\
Let $\bbbz/4\simeq H_{10}=\langle u\mid u^{4}=1\rangle$ and $\bbbz/4\simeq
H_{01}=\langle v\mid v^{4}=1\rangle$. We set $G=H_{10}\times H_{01}$, $H_{11}=\langle
u v\rangle$, $H_{21}=\langle u^{2} v\rangle$, $H_{12}=\langle u
v^{2}\rangle$, $H_{31}=H_{13}=\langle  v^{3}\rangle$, $T_{10}=\langle
u^{2}\rangle$, $T_{01}=\langle v^{2}\rangle$, $T_{11}=\langle u^{2}
v^{2}\rangle$. If $\chi_{1},\chi_{2}$ generate $({\bbbz}/4\times {\bbbz}/4)^{\star}$,
the building data of the cover must satisfy:
\begin{equation}
\label{zetaquattroperzetaquattro}
\left\{ \begin{array}{rl}
4{\cal{L}}_{\chi_{1}}=
&[\sum_{i=0}^{3}(D_{H_{1i}\phi_{1i}}+3D_{H_{1i}\phi_{1i}^{3}})]+
2(D_{H_{21}\phi_{21}}+ D_{H_{21}\phi_{21}^{3}}+D_{T_{10}}+ D_{T_{11}})\\
4{\cal{L}}_{\chi_{2}}=
&[\sum_{i=0}^{3}(D_{H_{i1}\phi_{i1}}+3D_{H_{i1}\phi_{i1}^{3}})]+
2(D_{H_{12}\phi_{12}}+ D_{H_{12}\phi_{12}^{3}}+D_{T_{01}}+ D_{T_{11}})\\
\end{array}\right.
\end{equation}
\noindent
By table (\ref{tabellauno}) and by \ref{otto} we can assume that the action on
$D$ has $V_{2,\chi_{2}^{3}}=0$, i.e. $l_{2}={\rm{deg}}L_{\chi_{2}}=1$. Now
we face two cases:
$V_{2,\chi_{1}^{3}}=0$ or $V_{2,\chi_{1}^{2}\chi_{2}}=0$. The case
$V_{2,\chi_{1}^{2}\chi_{2}}=0$ i.e. 
$l_{\chi_{1}^{2}\chi_{2}^{3}}={\rm{deg}}L_{\chi_{1}^{2}\chi_{2}^{3}}=1$ does not
occur. If $V_{2,\chi_{1}^{3}}=0$ that is 
$l_{\chi_{1}}={\rm{deg}}L_{\chi_{1}}=1$ then $p_{g}(S)\leq 2$.
\\
\\
{\bf{The case}} $G=\bbbz/2\times\bbbz/8$:\\
Let $\bbbz/2\simeq H=\langle u\mid u^{2}=1\rangle$ and $\bbbz/8\simeq
K=\langle v\mid v^{8}=1\rangle$. We set $G=H\times K$, $H_{2}=\langle
uv^{4}\rangle$, $K_{2}=\langle v^{4} \rangle$, $H_{4}=\langle u
v^{2}\rangle$, $K_{4}=\langle  v^{2}\rangle$, $H_{8}=\langle
uv\rangle$. If $\chi_{1},\chi_{2}$ generate $({\bbbz}/2\times {\bbbz}/8)^{\star}$,
the building data of the cover must satisfy:
\begin{equation}
\label{zetadueperzetaotto}
\left\{ \begin{array}{l}
2{\cal{L}}_{\chi_{1}}=D_{H_{}}+D_{H_{2}}+D_{H_{4}\phi_{1}}+D_{H_{4}\phi_{1}^{3}}
+\sum_{\!\!\!\!\!\!\!\!\!{_{_{_{_{_{_{\scriptscriptstyle{i=1,3,5,7}^{}}}}}}}}}
\!\!\!\!\!D_{H_{8}\phi^{i}}\\
\\

8{\cal{L}}_{\chi_{2}}=4(D_{H_{2}}\!\!\!+\!\!D_{K_{2}})\!\!+\!\!
2(D_{H_{4}\phi_{1}}\!\!\!+\!\!K_{H_{4}\psi_{1}})\!\!+\!6

(D_{H_{4}\phi_{1}^{3}}\!\!+\!\!D_{K_{4}\psi_{1}^{3}}\!)

+\sum{\!\!\!\!\!\!\!\!\!{_{_{_{_{_{_{\scriptscriptstyle{i=1,3,5,7}^{}}}}}}}}}
\!\!\!\!\!i(D_{H_{8}\phi^{i}}\!\!+\!\!D_{K\psi^{i}}\!)\\
\end{array}\right.
\end{equation}
\noindent
We know that, up to $G$ isomorphisms, the action on $F$ has three branch
points on $A$ with the following action: 
$d^{1}_{H_{2}}={\rm{deg}}D_{H_{2}}=1$,
$d^{1}_{K\psi^{7}}={\rm{deg}}D_{K\psi^{7}}=1$,
$d^{1}_{H_{8}\phi^{5}}={\rm{deg}}D_{H_{8}\phi^{5}}=1$. To compute the suitable
actions on $D$ we can assume $V_{2,\chi_{2}^{3}}=0$ that
is ${\rm{deg}}{\cal{L}}_{\chi_{2}}=1$. Now we have to discuss two
cases: $V_{2,(\chi_{1}\chi_{2}^{2})^{-1}}=0$ or $V_{2,\chi_{2}^{-3}}=0$.\\

If $V_{2,(\chi_{1}\chi_{2}^{2})^{-1}}=0$
then ${\rm{deg}}{\cal{L}}_{\chi_{1}\chi_{2}^{2}}=1$ and with this condition the system
(\ref{zetadueperzetaotto}) does not have any solution.\\

If $V_{2,\chi_{2}^{-3}}=0$ then we have to add to (\ref{zetadueperzetaotto}) the
condition ${\rm{deg}}{\cal{L}}_{\chi_{2}^{3}}=1$. We set
${\rm{deg}}{\cal{L}}_{\chi_{1}}=m$ and we find the following solutions: 

\begin{equation}
\label{soluzionisistemauno}
\begin{tabular}{|l|l|}\hline
  $i) $ &  ${\rm{Decomposition}\, of }\,H^{0}(D,\Omega^{1}_{D})$\\
\hline
$d^{2}_{K\psi^{7}}=d^{2}_{K\psi}=d^{2}_{H_{8}\phi^{7}}=d^{2}_{H_{8}\phi}=0$ &
$l_{\chi_{2}^{i}}=1$, $i=1,...,7$
\\ \hline
$
d^{2}_{H_{4}\phi_{1}}=d^{2}_{H_{4}\phi_{1}^{3}}=
d^{2}_{K_{4}\psi_{1}^{}}=d^{2}_{K_{4}\psi_{1}^{3}}=0$
&
$l_{\chi_{1}\chi_{2}^{}}=l_{\chi_{1}\chi_{2}^{3}}=
l_{\chi_{1}\chi_{2}^{6}}=m+1-d^{2}_{H_{8}\phi^{5}}$\\
\hline

$d^{2}_{H_{2}}=d^{2}_{K_{2}}=0$ & 
$l_{\chi_{1}\chi_{2}^{2}}=l_{\chi_{1}\chi_{2}^{5}}=
l_{\chi_{1}\chi_{2}^{7}}=m+1-d^{2}_{H_{8}\phi^{3}}$ \\ \hline

$d^{2}_{K\psi^{3}}+d^{2}_{H_{8}\phi^{3}}=1$
& $l_{\chi_{1}\chi_{2}^{4}}=m+1-d^{2}_{H_{8}\phi^{3}}-d^{2}_{H_{8}\phi^{5}}$ \\\hline

$d^{2}_{K\psi^{5}}+d^{2}_{H_{8}\phi^{5}}=1$
&$g(D)=8m-1-4(d^{2}_{H_{8}\phi^{3}}+d^{2}_{H_{8}\phi^{5}})$\\\hline
\end{tabular}
\end{equation}
\noindent
\begin{equation}
\label{soluzionisistemadue}
\begin{tabular}{|l|l|}\hline
  $ii) $ &  ${\rm{Decomposition}\, of }\,H^{0}(D,\Omega^{1}_{D})$\\
\hline
$d^{2}_{K\psi^{i}}=d^{2}_{H_{8}\phi^{i}}=0$ $i=1,5,7$&
$l_{\chi_{2}^{i}}=1$, $i=1,3,4,6$, $l_{\chi_{2}^{j}}=2$, $j=2,5,7$
\\ \hline

$
d^{2}_{H_{4}\phi_{1}^{3}}=d^{2}_{K_{4}\psi_{1}^{3}}=0$
&
$l_{\chi_{1}\chi_{2}^{2}}=l_{\chi_{1}\chi_{2}^{7}}=m+2-d^{2}_{H_{8}\phi^{3}}-
d^{2}_{H_{4}\phi_{1}}$\\
\hline

$d^{2}_{H_{2}}=d^{2}_{K_{2}}=0$ & 
$l_{\chi_{1}\chi_{2}^{3}}=l_{\chi_{1}\chi_{2}^{6}}=m+1-d^{2}_{H_{4}\phi_{1}}$,
$l_{\chi_{1}\chi_{2}}=m+1$
\\
\hline

$d^{2}_{K\psi^{3}}+d^{2}_{H_{8}\phi^{3}}=2$
&
$l_{\chi_{1}\chi_{2}^{5}}=l_{\chi_{1}\chi_{2}^{4}}=
m+1-d^{2}_{H_{8}\phi^{3}}$
\\\hline

$d^{2}_{K_{4}\psi_{1}}+d^{2}_{H_{4}\phi_{1}}=1$
&$g(D)=8m+4-4(d^{2}_{H_{8}\phi^{3}}+d^{2}_{H_{4}\phi_{1}})$\\\hline
\end{tabular}
\end{equation}
\noindent

\begin{equation}
\label{soluzionisistematre}
\begin{tabular}{|l|l|}\hline
  $iii) $ &  ${\rm{Decomposition}\, of }\,H^{0}(D,\Omega^{1}_{D})$\\
\hline
$d^{2}_{K\psi^{i}}=d^{2}_{H_{8}\phi^{i}}=0$ $i=5,7$ &
$l_{\chi_{2}^{i}}=1$, $i=1,2,3,4,6$, $l_{\chi_{2}^{j}}=2$, $j=5,7$,
\\ \hline

$
d^{2}_{H_{4}\phi_{1}^{i}}=d^{2}_{K_{4}\psi_{1}^{i}}=0$, $i=1,3$
&
$l_{\chi_{1}\chi_{2}^{}}=l_{\chi_{1}\chi_{2}^{3}}=m+1-d^{2}_{H_{2}}$\\
\hline

$d^{2}_{H_{2}}+d^{2}_{K_{2}}=1$ & 
$l_{\chi_{1}\chi_{2}^{2}}=l_{\chi_{1}\chi_{2}^{4}}=m+1-d^{2}_{H_{8}\phi^{3}}$,
$l_{\chi_{1}\chi_{2}^{6}}=m+1$ \\
\hline

$d^{2}_{K\psi^{3}}+d^{2}_{H_{8}\phi^{3}}=1$
&
$l_{\chi_{1}\chi_{2}^{5}}=l_{\chi_{1}\chi_{2}^{7}}=
m+2-d^{2}_{H_{8}\phi^{3}}-d^{2}_{H_{2}}$
\\\hline

$d^{2}_{K\psi}+d^{2}_{H_{8}\phi}=1$
&$g(D)=8m+3-4(d^{2}_{H_{8}\phi^{3}}+d^{2}_{H_{2}})$\\\hline
\end{tabular}
\end{equation}
\noindent

\begin{equation}
\label{soluzionisistemaquattro}
\begin{tabular}{|l|l|}\hline
  $iv) $ &  ${\rm{Decomposition}\, of }\,H^{0}(D,\Omega^{1}_{D})$\\
\hline
$d^{2}_{K\psi^{i}}=d^{2}_{H_{8}\phi^{i}}=0$ $i=5,7$ &
$l_{\chi_{2}^{1}}=l_{\chi_{2}^{3}}1$,
$l_{\chi_{2}^{2}}=l_{\chi_{2}^{4}}=l_{\chi_{2}^{6}}=2$,$l_{\chi_{2}^{5}}=
l_{\chi_{2}^{7}}=3$,
\\ \hline

$
d^{2}_{H_{4}\phi_{1}^{i}}=d^{2}_{K_{4}\psi_{1}^{i}}=0$, $i=1,3$
&
$l_{\chi_{1}\chi_{2}^{}}=l_{\chi_{1}\chi_{2}^{3}}=m+1$\\
\hline

$d^{2}_{H_{2}}=d^{2}_{K_{2}}=0$ & 
$l_{\chi_{1}\chi_{2}^{2}}=l_{\chi_{1}\chi_{2}^{4}}=m+2-d^{2}_{H_{8}\phi^{3}}$,
$l_{\chi_{1}\chi_{2}^{6}}=m+2$ \\
\hline

$d^{2}_{K\psi^{3}}+d^{2}_{H_{8}\phi^{3}}=2$
&
$l_{\chi_{1}\chi_{2}^{5}}=l_{\chi_{1}\chi_{2}^{7}}=
m+3-d^{2}_{H_{8}\phi^{3}}$
\\\hline

$d^{2}_{K\psi}+d^{2}_{H_{8}\phi}=2$
&$g(D)=8m+12-4d^{2}_{H_{8}\phi^{3}}$\\\hline
\end{tabular}
\end{equation}
\noindent It is an easy task to show that they are the claimed ones. Moreover to
obtain the claimed result in table (\ref{mostro}) we have to analize the
singularities; a straightforward long computation. We remark only that in table
(\ref{mostro}) appears some entries with the same value for $p_{g}$ and $g(D)$, but
different $K^{2}_{S}$. In fact this depends only on the type of the singular points. 
\qed\\

\begin{theorem}\label{eccolottouno}
If $S$ is a $G$-sandwich canonically $3$-fibred surface on an elliptic curve
with $\mid G\mid>8$, $p_{g}(S)=m\geq 2$ then $q(S)=1$, $S=Z/G$ where
$G={\bbbz}/2\times{\bbbz}/8$ acts diagonally and without fixed points on
$Z=F\times D$, $g(F)=3$, $g(D)=8m+1$.
\end{theorem}
\proof The same argument as for Theorem \ref{eccolotto}, but much more easier.\qed\\
We have shown the second theorem stated in the introduction; we will explain how to
look these new families.

\subsubsection{A geometrical construction of the family of \ref{eccolottouno}} 
The families in theorem \ref{eccolotto} have a rich geometry which should be
studied. Here we will
show the easier geometry of the family in
\ref{eccolottouno}. We must understand first the
${\bbbz}/2\times{\bbbz}/8$-action on the genus-$3$ curve
$F$ whose existence is stated in \ref{diciannove}. We fix three distinct
points
$Q_{1}, Q_{2}, Q_{3}\in A=I\!\!P^{1}$ and we consider the $2$-to-$1$ cover
branched on $Q_{2}, Q_{3}$, $\rho_{\Gamma_{1}}:\Gamma_{1}\rightarrow A$ and
the $8$-to-$1$ cover branched on $Q_{1},Q_{2},Q_{3}$,
$\rho_{\Gamma_{2}}:\Gamma_{2}\rightarrow A$ such that: 
$\rho_{\Gamma_{1}}^{-1}(Q_{2})=Q_{2}^{1}$, 
$\rho_{\Gamma_{1}}^{-1}(Q_{3})=Q_{3}^{1}$, 
$\rho_{\Gamma_{1}}^{-1}(Q_{1})=\{P_{1}^{1}, P_{2}^{1}\}$, 
$\rho_{\Gamma_{2}}^{-1}(Q_{1})=Q_{1}^{2}$,
$\rho_{\Gamma_{2}}^{-1}(Q_{3})=Q_{3}^{2}$, 
$\rho_{\Gamma_{2}}^{-1}(Q_{2})=\{P_{1}^{2}, P_{2}^{2},P_{3}^{2}, P_{4}^{2}\}$. If  
$i_{\Gamma_{1}}:\Gamma_{1}\rightarrow \Gamma_{1} $ and 
$i_{\Gamma_{2}}:\Gamma_{2}\rightarrow \Gamma_{2} $ and are the generators of the
deck-transformation groups of   $\rho_{\Gamma_{1}}$ and $\rho_{\Gamma_{2}}$
respectively, then the group $\langle i_{\Gamma_{1}}\rangle\times\langle
i_{\Gamma_{2}}\rangle\simeq{\bbbz}/2\times{\bbbz}/8$ acts on the
normalization $F$ of the fibre product
$\Gamma_{1}\times_{A}\Gamma_{2}$,
$F/\langle i_{\Gamma_{1}}\rangle=\Gamma_{2}$,
$F/\langle i_{\Gamma_{2}}\rangle=\Gamma_{1}$. We denote by
$\pi_{F}:F\rightarrow A$ the Galois map and we
set:
$\pi_{F}^{-1}(Q_{1})=\{A_{1},A_{2}\}$,
$\pi_{F}^{-1}(Q_{3})=\{B_{1},B_{2}\}$,
$\pi_{F}^{-1}(Q_{2})=\{S_{1},...,S_{8}\}$. Then $G(A_{1})=G(A_{2})=\langle
i_{\Gamma_{2}}\rangle$, $G(B_{1})=G(B_{2})=\langle
i_{\Gamma_{1}}i_{\Gamma_{2}}\rangle$, and $G(S_{i})=\langle
i_{\Gamma_{1}}i_{\Gamma_{2}}^{4}\rangle$, $i=1,...,8$. It is easy to see that the
three canonical divisors:
${\rm{div}}\eta_{1}=2A_{1}+2A_{2}$,
${\rm{div}}\eta_{2}=2B_{1}+2B_{2}$, ${\rm{div}}\eta_{3}=A_{1}+A_{2}+B_{1}+B_{2}$,
give a basis $\langle\eta_{1},\eta_{2},\eta_{3}\rangle$ of
$H^{0}(F,\Omega^{1}_{F})$ such
that:
$\langle\eta_{1}\rangle=V_{1,\chi_{2}}$,
$\langle\eta_{2}\rangle=V_{1,\chi_{2}^{3}}$,
$\langle\eta_{3}\rangle=V_{1,\chi_{1}\chi_{2}^{2}}$. On the other hand the $G$-action
on $D$ with elliptic quotient
$B$ is easily obtained from the fibre product
${\widetilde{\Gamma_{1}}}\times_{B}{\widetilde{\Gamma_{2}}}$
where ${\widetilde{\Gamma_{1}}}\rightarrow B$ is the $2$-to-$1$ cover branched on
$2m$ points $R_{1},...,R_{2m}$ while ${\widetilde{\Gamma_{2}}}\rightarrow B$ is
an unramified cover of degree $8$. In fact it is easy to obtain a
basis $\langle\epsilon\rangle\cup\langle\alpha_{i}^{j}\rangle_{i=1,...,m}^{j=0,...,7}$
of $H^{0}(F,\Omega^{1}_{F})$ such
that:
$i_{{\widetilde{\Gamma_{1}}}}^{\star}\epsilon=\epsilon$,
$i_{{\widetilde{\Gamma_{1}}}}^{\star}\alpha_{i}^{j}=-\alpha_{i}^{j}$, $i=1,...,m,
j=0,...,7$,
$i_{{\widetilde{\Gamma_{2}}}}^{\star}\epsilon=\epsilon$,
$i_{{\widetilde{\Gamma_{2}}}}^{\star}\alpha_{i}^{j}=\mu^{j}\alpha_{i}^{j}$,
$i=1,...,m, j=0,...,7$ where $\mu$ is a primitive $8$-root of the unity. In particular
$V_{2{\rm{id}}}= \langle\epsilon \rangle$, $V_{2\chi_{2}^{j}}=0$, $j=0,...,7$
and $V_{2\chi_{1}\chi_{2}^{j}}=\langle \alpha_{i}^{8-j}\rangle_{j=0,...,7}$,
$i=1,...,m$. Now it is trivial to see that the ${\bbbz}/2\times{\bbbz}/8$-diagonal
action induced by the group $G=\langle(i_{\Gamma_{1}}i_{{\widetilde{\Gamma_{1}}}}),
(i_{\Gamma_{2}}i_{{\widetilde{\Gamma_{2}}}})\rangle$ on the product surface 
$Z=F\times D$ gives a quotient surface $S=Z/G$ such
that: $H^{0}(X,\Omega^{1}_{X})=H^{0}(Z,\Omega^{2}_{Z})^{G}\simeq
\langle\epsilon\rangle$ and $H^{0}(X,\Omega^{2}_{X})=V_{1,\chi_{1}\chi_{2}^{2}}\otimes
V_{2,\chi_{1}\chi_{2}^{6}}=\{\eta_{3}\wedge\alpha_{1}^{6},...,\eta_{3}
\wedge\alpha_{m}^{6}\}$. In particular $m=p_{g}(S)$ and the image of the canonical map
$\Phi_{\mid K_{S}\mid}$ is the normal elliptic curve in $I\!\!P^{p_{g}-1}$ of degree
$p_{g}$.\qed\\\\
We end this article with a natural problem: {\it{classify all the isotrivial canonical
fibrations where the involved group $G$ is not Abelian}}. 


{\samepage}
Zucconi Francesco\\
Universit\`{a} di Udine Dipartimento di Matematica e Informatica\\
Via delle Scienze 206 33100 Udine, Italia\\
e-mail zucconi@dimi.uniud.it
\end{document}